\documentclass[10pt,leqeq]{article}
\usepackage{amssymb,amsthm}
\usepackage{amsmath}
\usepackage{amscd}
\usepackage{epsf,graphicx}
\usepackage[all]{xy}
\newtheorem{thm}{Theorem}
\newtheorem{lema}[thm]{Lemma}
\newtheorem{cor}[thm]{Corollary}
\newtheorem{prop}[thm]{Proposition} 

\newtheorem{defi}[thm]{Definition}
\newtheorem{exa}[thm]{Example}

\newcommand{\des}{\displaystyle}

\renewcommand{\ker}{{\rm{Ker}}}
\newcommand{\en}{{\rm{End}}}
\newcommand{\id}{{\rm{id}}}
\newcommand{\vol}{{\rm{vol}}}
\newcommand{\sgn}{{\rm{sgn}}}
\newcommand{\sym}{{\rm{Sym}}}
\newcommand{\qsym}{{\rm{QSym}}}

\newcommand{\sdim}{\rm{sdim}}
\newcommand{\inv}{{\rm{Inv}}}
\newcommand{\coinv}{{\rm{Coinv}}}

\newcommand{\ocat}{{\bf{OCat}}}
\newcommand{\oalg}{{\bf{OAlg}}}
\newcommand{\pcat}{{\bf{Pre\mbox{-}Cat}}}
\newcommand{\cat}{{\bf{Cat}}}
\newcommand{\fu}{{\rm{Funct}}}
\newcommand{\func}{{\rm{Func}}}
\newcommand{\schur}{{\rm{Schur}}}

\newcommand{\mat}{{\rm{Mat}}}
\newcommand{\mo}{{\rm{Mor}}}

\newcommand{\alg}{{\rm{alg}}}
\newcommand{\ob}{{\rm{Ob}}}

\newcommand{\dif}{{\rm{d}}}

\newcommand{\cH}{{{\rm{H}}^{\bullet}}}
\newcommand{\Ho}{{{\rm{H}}_{\bullet}}}
\newcommand{\aut}{{\rm{Aut}}}
\newcommand{\vect}{{\bf{Vect}}}
\newcommand{\ho}{{\rm{Hom}}}
\renewcommand{\des}{\displaystyle}

\begin{document}
\setlength{\baselineskip}{16pt}
\title{Quantum symmetric functions}
\author{Rafael D\'\i az\thanks{Work partially supported by UCV.}\ \
and Eddy Pariguan\thanks{Work partially supported by FONACIT.}}
\maketitle
\begin{abstract}
We study quantum deformations of Poisson orbivarieties. Given a
Poisson manifold $(\mathbb{R}^{m},\alpha)$ we consider the Poisson
orbivariety $(\mathbb{R}^{m})^{n}/S_{n}$. The Kontsevich star
product on functions on $(\mathbb{R}^{m})^{n}$ induces a star
product on functions on $(\mathbb{R}^{m})^{n}/S_{n}$. We provide
explicit formulae for the case ${{\mathfrak h} \times {\mathfrak
h}}/\mathcal{W}$, where ${\mathfrak h}$ is the Cartan subalgebra
of a classical Lie algebra ${\mathfrak g}$ and $\mathcal{W}$ is
the Weyl group of ${\mathfrak h}$. We approach our problem from a
fairly general point of view, introducing Polya functors for
categories over non-symmetric Hopf operads.
\end{abstract}

\section{Introduction}
Let $k$ be a field of characteristic $0$, $M$ be a set and $G$ be
a subgroup of the permutation group on $n$-letters $S_n$. A
function $f:M^n \to k$ is said to be $G$-symmetric if
$f(x_{\sigma(1)},x_{\sigma(2)},...,x_{\sigma(n)})=
f(x_1,x_2,...,x_n)$ for all $\sigma \in G\subseteq S_n$ and all
$x_1,...,x_n \in M$. The $k$-space of $G$-symmetric functions
${\func}(M^n,{k})^{G}$ is a subalgebra of the $k$-algebra
${\func}(M^n,k)$ of all functions from $M^n$ to $k$. One of the
goals of this paper is to find explicit formulae for a product on
the algebra ${\func}(M^n,k)^{G}$ in a variety of contexts. Our
approach is based on the following observations:
\begin{itemize}
\item{ It is often easier to work with coinvariant functions
${\func}(M^n,k)_G$ instead of working with invariant functions.}
\item{ Symmetric functions arise as an instance of a general
construction which assigns to any $k$-algebra $A$ its $n$-th
symmetric power algebra $\sym^{n}(A)$. This insight led us to
introduce the notion of Polya functors which we present in the
context of categories over non-symmetric Hopf operads.}
\end{itemize}
Our main interest is to study formal deformations of the algebra
${\func}(M^n,k)^{G}$. We take the real numbers $\mathbb{R}$ as the
ground field, and let $(\mathbb{R}^{m},\{\mbox{},\mbox{}\})$ be a
Poisson manifold. Under this conditions Kontsevich in \cite{Kon}
have shown the existence of a canonical formal deformation
$(C^{\infty}(\mathbb{R}^{m})[[\hbar]],\star )$ of the algebra
$(C^{\infty}(\mathbb{R}^{m}),\cdot )$ of smooth functions on
$\mathbb{R}^{m}$. We prove that if the Poisson bracket on
$(\mathbb{R}^{m},\{\mbox{},\mbox{}\})$ is $G$-equivariant for
$G\subset S_m$, then the $\star$-product on
$C^{\infty}(\mathbb{R}^{m})[[\hbar]]$ induces a $\star$-product on
the algebra of symmetric functions
$C^{\infty}(\mathbb{R}^{m})^{G}[[\hbar]]$, which we call the
algebra of {\it quantum symmetric functions}. We regard this
algebra as the deformation quantization of the Poisson orbifold
$\mathbb{R}^{m}/G$. We remark that in a recent paper \cite{Dol},
Dolgushev has proved the existence of a quantum product on the
algebra of invariant functions $C^{\infty}(M)^{G}[[\hbar]]$ for an
arbitrary Poisson manifold $M$ acted upon by a finite group $G$.
His result is based on an alternative proof of the Kontsevich
formality theorem which is manifestly  covariant.

We present a general description of the quantum product on
$\mathbb{R}^{m}/G$ using the Kontsevich star product. We give
explicit formulae for the product rule in the following three
cases:
\begin{itemize}
\item{ symplectic orbifold ${\mathfrak h}\times {\mathfrak
h}/{\mathcal{ W}}$ where ${\mathfrak h}$ is a Cartan subalgebra of
a classical Lie algebra ${\mathfrak g}$, and $\mathcal{W}$ is the
Weyl group associated to ${\mathfrak h}$,} \item{ symplectic
orbifold $\mathbb{C}^{n}/{\mathbb{Z}_m^{n}\rtimes S_n}$, } \item{
symplectic orbifold $\mathbb{C}^{n}/{\mathcal{D}_m^{n}\rtimes
S_n}$, where $\mathcal{D}_m$ is the dihedral group of $2m$
elements.}
\end{itemize}

Our motivation to consider these orbifolds came from the study of
noncommutative solitons in orbifolds $\cite{Min2}$, $\cite{EM}$
and the quantization of the moduli space of vacua in $M$-theory as
consider in the matrix model approach. Our results will be raised,
to the categorical context in $\cite{DP2}$. In a different
direction, they may be extended to include the $q$-Weyl and
$h$-Weyl algebras as it is done in $\cite{DP1}$. We would like to
mention that these orbifolds have also been studied from a
different point of view in $\cite{Al}$, and more recently in
$\cite{Etg2}$.

We consider the quantum symmetric functions of type $A_n$ and
uncover its relation with the $\schur(\infty,n)$ algebras. The
latter algebras are natural generalizations of the Schur algebras
as defined in $\cite{JG}$. We also study the symmetric powers of
the $M$-Weyl algebra, which we define as the algebra generated by
$x^{-1}$ and $\frac{\partial}{\partial x}$. We provided explicit
formulae for the normal coordinates for the $M$-Weyl algebra as
well as for its symmetric powers.  Similarly, we make clear the
relation between quantum symmetric odd-functions and the Schur
algebras of various dimensions . Finally, we give a cohomological
interpretation of the algebra of supersymmetric functions.

\section{Invariants vs coinvariants}\label{sinco}

In this section we introduce the notion of Polya functors for
categories over non-symmetric Hopf operads, and provide a list of
applications of the Polya functors. We will consider invariant
theory for finite groups as well as for compact topological
groups. To avoid duplication we will consider only the
latter case in the proofs.\\
Let $k$ be a field of characteristic $0$ and consider
$(\vect_k,\otimes,k)$ the monoidal category of vector spaces with
linear transformations as morphisms. For any set $I$, consider the
category of $I$-graded  vector spaces $\vect_I$, it has as objects
$I$-graded vector spaces, ${\des V=\bigoplus_{i\in I} V_i}, \ \ \
V_i\in \ob(\vect_k).$  Morphisms between objects $V,W\in
\ob(\vect_I)$ are given by $\mo (V,W)={\des\prod_{i\in I}\ho
(V_i,W_i)}.$ The category $(\vect_I, \otimes_I,k_I)$ has a
monoidal structure compatible with direct sums induced by the
corresponding  structures on $(\vect_k,\otimes, k)$. Explicitly,
given $V,W\in \ob(\vect_I)$ we have $(V\oplus W)_i=V_i\oplus W_i$,
$(V\otimes W)_i=V_i\otimes W_i$, and $(k_I)_i=k$. For a finite
group $G$, we denote by $\vect_I(G)$ the category of $I$-graded
vector spaces provided with grading preserving $G$ actions.
Morphisms in $\vect_I(G)$ are intertwiners, i.e., maps $\varphi:
V\longrightarrow W$ such that $\varphi(gv)=g\varphi(v)$, for all
$v\in V, g\in G$. Abusing notation, for an infinite compact
topological group provided with a biinvariant Haar measure $\dif
g$, we denote by $\vect_I(G)$ the category of finite dimensional
vector spaces over $\mathbb{C}$ provided with a $G$ action. We
define the symmetrization map $s_V:V \longrightarrow V^{G}$ as the
map given by
\begin{eqnarray*}
s_V(v)&=&\frac{1}{\mbox{vol}(G)}\int_{G}(gv)\dif g,\ \ \ \
\mbox{if}\ \ G \ \ \mbox{is infinite and}\ \ k=\mathbb{C},\\
\mbox{}& \mbox{}& \mbox{}\\
s_V(v)&=& \frac{1}{\sharp(G)}\sum_{g\in G}gv,\ \ \ \mbox{if}\ \ G
\ \ \mbox{is finite and}\ \  k \ \ \mbox{is a field of
characteristic zero,}
\end{eqnarray*}\\
where $\sharp(G)$ denotes the cardinality of $G$. The sequence  $0
\longrightarrow
\ker(s_V)
\longrightarrow V\longrightarrow V^{G}\longrightarrow 0,$ \ is
exact and  we obtain the corresponding commutative triangle
\[
\xymatrix{ V \ar[rr] \ar[dr] & & V^G
\\ & V/\ker(s_V) \ar[ur]_{s_V} & }
\]
\\
We denote the space $V/\ker(s_V)$ by $V_G$, and thus we have an
isomorphism $s_V:V_{G} \longrightarrow V^{G}$. We define two
functors
$$\begin{array}{cccc}
  {\rm{Inv}}: & \vect_I(G) & \longrightarrow & \vect_I \\
  \mbox{ } & V & \longmapsto & V^{G}
\end{array}  $$
$$\begin{array}{cccc}
  {\rm{Coinv}}: & \vect_I(G) & \longrightarrow & \vect_I \\
  \mbox{ } & V & \longmapsto & V_{G}
\end{array}  $$
For any $v\in V$, $\overline{v}$ denotes the equivalence class of
$v$ in $V_G$. We have the following
\begin{prop}\label{inco}
The maps $s_V$ above define a natural isomorphism $s:{\rm
{Coinv}\longrightarrow {\rm{Inv}}}$.
\end{prop}

\begin{proof} For each $V\in  \vect_I(G)$ the construction above
provides an isomorphism $s_V: V_G\longrightarrow V^G$. For a given
morphism $V \xrightarrow {\alpha} W$, we have
\begin{eqnarray*}
\alpha \circ s_V(\overline{v})&=& \alpha\left(
\frac{1}{\vol(G)}\int_G(gv)\dif
g\right)=\frac{1}{\vol(G)}\int_G \alpha(gv)\dif g\\
\mbox{}& \mbox{} &\mbox{} \\
 &=&\frac{1}{\vol(G)}\int_G g(\alpha v)\dif g=s_W (\overline{\alpha v})=s_W
 \overline{\alpha}(\overline{v}),\ \ \mbox{for all}\ \ v\in V
\end{eqnarray*}
thus proving that for each arrow $V \xrightarrow {\alpha} W$ the
diagram
\[
\xymatrix{ V_G \ar[r]^{s_V} \ar[d]_{{\rm {Coinv(\alpha)}}} & V^G
\ar[d]^{{\rm {Inv(\alpha)}}}
\\ W_G \ar[r]^{s_W} & W^G }
\]
is commutative.\end{proof}
Notice that $V_G$ may  also be defined
as $V_{G}=V/\langle v-gv: \
\ v\in V,\  g\in G\rangle$.  Constructions above can be
generalized to the category ${\bf{Cat}}_k$ of all $k$-linear
categories. We make a more general construction in the next
section in order to include linear categories over non-symmetric
Hopf operads.

\subsection{Categories over non-symmetric operads}

We review the notion of operads and the notion of algebras over
operads $\cite{May}$. For finite groups $G\subset H$ such that $G$
acts on a $k$-vector space $V$, the induced representation is
defined by ${\rm {Ind}}_G^{H}(V)=(k[G]\otimes V)_H$ where $k[G]$
denotes the group algebra of $G$. Let us define two $k$-linear
categories $\mathbb{N}$ and $S$

\begin{center}
\begin{tabular}{lll}
$ \bullet \ \ob(\mathbb{N})=\{0,1,2,\dots,n,\dots\}$ & \hspace{0.5cm} & $\bullet \ \ob(S)=\{0,1,2,\dots,n,\dots\}$\\
\mbox{} & \mbox{}\\
$\bullet \ \mo_{\mathbb{N}}(n,m)=\left \{ \begin{array}{cc}
  k, & \mbox{if}\ m=n \\
  0, & \mbox{if}\ m\neq n \\
\end{array}\right.$ & \hspace{0.5cm}& $\bullet \ \mo_{S}(n,m)=\left \{ \begin{array}{cc}
S_n, & \mbox{if}\ m=n \\
  0, & \mbox{if}\ m\neq n \\
\end{array}\right.$
\end{tabular}
\end{center}

The category $\fu(S^{op},\vect_k)$ of contravariant functors from
$S$ to $\vect_k$ possesses three important monoidal structures
given on objects by
\begin{itemize}
\item{$(V+W)(n)=V(n)\oplus W(n)$.} \item{$V\otimes W(n)={\des
\bigoplus_{i+j=n} {\rm {Ind}} _{S_i\times S_j} ^{S_n}
(V(i)\otimes_k W(j))}$.} \item{$V \circ W(n)={\des
\bigoplus_{p\geq 0}\bigoplus_{a_1+\dots+a_p=n}{\rm
{Ind}}_{S_{a_1}\times \dots \times S_{a_p}} ^{S_n}
V(p)\otimes_{s_p} W(a_1)\otimes_{k}\dots \otimes_k W(a_p)}$.}
\end{itemize}
 The category $\fu(\mathbb{N},\vect_k)$ admits similar monoidal
structures by forgetting the $S_n$ actions.

\begin{defi}
\begin{itemize}
\item{ An operad is a monoid in the monoidal category
$(\fu(S^{op},\vect_k),\circ,1)$, where $1(n)=0, n\neq 1$ and
$1(1)=k$.} \item{ A nonsymmetric operad is a monoid in the
monoidal category $(\fu(\mathbb{N},\vect_k),\circ,1)$.}
\end{itemize}
\end{defi}

Explicitly, an (non-symmetric) operad is given by $m_p:
\mathcal{O}(p)\otimes \mathcal{O}(a_1)\otimes \dots \otimes \mathcal{O}(a_p)\longrightarrow \mathcal{O}(a_1+\dots+a_p)$
satisfying the list of axioms given for example in $\cite{May}$.
If no confusion arises we write $m$ instead of $m_p$.
\begin{defi}
Let $V\in \vect_k$, we define the endomorphisms operad by
$\en_V(n)=\ho(V^{\otimes n}, V)$, for all $n\in\mathbb{N}$.
Composition are given by arrows
\[
\xymatrix{\ho(V^{\otimes p},V)\otimes\ho(V^{\otimes a_1},V)\otimes
\dots\otimes \ho(V^{\otimes a_p},V) \ar[d] \\
\ho(V^{\otimes p},V)\otimes \ho(V^{\otimes (a_1+\dots+a_p)},V^{\otimes p})
\ar[d]\\
 \ho(V^{\otimes n},V)}
\]
for integers $n,a_1,\dots,a_p$ such that $a_1+\cdots+a_p=n$. For
more details see $\cite{May}$.

\end{defi}

Let us introduce the category $\pcat_k$ of small pre-categories.
Objects of $\pcat_k$ are called small pre-categories. A small
pre-category $\mathcal{C}$ consists of the following data:
\begin{itemize}
\item{A set of objects $\ob(\mathcal{C})$.}
\item{A vector space $\mo_{\mathcal{C}}(x,y)$ associated to each pair of objects $x,y\in\ob(\mathcal{C})$.}
\end{itemize}
A morphisms $F\in\mo_{\pcat}(\mathcal{C},\mathcal{D})$ from
pre-category $\mathcal{C}$ to pre-category $D$ consists of a map\\
$F:\ob(\mathcal{C})\longrightarrow \ob(\mathcal{D})$ and a family
of maps $F_{x,y}: \mo_{\mathcal{C}}(x,y)\longrightarrow
\mo_{\mathcal{D}}(F(x),F(y))$, for $x,y\in\ob(\mathcal{C})$.\\
$\pcat_k$ has a natural partial monoidal structure. Given
pre-categories $\mathcal{C}$ and $\mathcal{D}$ such that\\
$\ob(\mathcal{C})=\ob(\mathcal{D})=X$, we define the product
pre-category $\mathcal{C} \mathcal{D}$ as follows
\begin{itemize}
\item{$\ob(\mathcal{C} \mathcal{D})=X$.}
\item{$\mo_{\mathcal{C} \mathcal{D}}(x,z)=\bigoplus_{y\in X}\mo_{\mathcal{C}}(x,y)\otimes \mo_{\mathcal{D}}(y,z)$.}
\end{itemize}
The partial units are the pre-categories $k_X$, defined as follows
\begin{itemize}
\item{${\ob}(k_X)=X$.} \item{${\mo}_{k_X}(a,b)=\left\{
\begin{array}{cc}
  k, & \ \mbox{if}\ \ a=b \\
  0, & \mbox{otherwise}. \\
\end{array}\right. $}
\end{itemize}
Given a pre-category $\mathcal{C}$, we define the non-symmetric
operad
$\en_{\mathcal{C}}(n)=\mo_{\pcat}(\mathcal{C}^{n},\mathcal{C})$,
$n\in \mathbb{N}$. We used the convention
$\mathcal{C}^{0}=k_{\ob(\mathcal{C})}$.
\begin{defi}\label{def5}
Let $\mathcal{O}$ be a non-symmetric $k$-linear operad. An
$\mathcal{O}$-category $(\mathcal{C},\gamma)$ is a pre-category
$\mathcal{C}$ together with a non-symmetric operad morphism
$\gamma:\mathcal{O}\longrightarrow \en_{\mathcal{C}}$. Explicitly,
a $k$-linear $\mathcal{O}$-category $\mathcal{C}$ consist of the
following data:
\begin{itemize}
\item{Objects $\ob(\mathcal{C})$.} \item{Morphisms:
$\ho_\mathcal{C}(x,y)\in \ob(\vect_k)$,  for each pair $x,y\in
\ob(\mathcal{C})$.} \item{For each, $k\in \mathbb{N}$ and objects $x_0,x_1,\dots,x_k\in
\ob(\mathcal{C})$ maps
$$p_{x_0,\dots,x_k}:{\mathcal{O}(k)}\otimes\ho_\mathcal{C}(x_0,x_1)\otimes \dots \otimes \ho_\mathcal{C}(x_{k-1},x_k)
\longrightarrow \ho_\mathcal{C}(x_0,x_k)$$ We usually write $p$
instead of $p_{x_0,\dots,x_k}$.}
\end{itemize}
These data should satisfy the following associativity axiom: Given
objects $x_0,\dots,x_{n_1+\dots+n_k}$, and morphisms
$a_i\in\ho_\mathcal{C}(x_{i-1},x_i), i\in[n_1+\dots+n_k]$, $t\in
{\mathcal{O}(k)}$, $t_i\in {\mathcal{O}(n_i)},$ then
$$p(m(t;t_1, \dots, t_k); a_{1},  \dots, a_{n_1+\dots+n_{k}})=$$
$$p(t, p(t_1; a_1,\dots, a_{n_1}),\dots, p(t_k; a_{n_{1}+\dots+
n_{(k-1)}+1},\dots, a_{n_1+\dots +n_k}) ).$$
\end{defi}
\newpage

For example if we are given objects $x_i\in\ob(\mathcal{C})$, for
$i=0,1,\dots,5$, morphisms  $a_i\in\mo(x_{i-1},x_i)$ for
$i=1,\dots,5$ and $t\in\mathcal{O}(5)$, then the morphism
$p(t;a_1,\dots,a_5)$ from object $x_0$ to object $x_5$ is
represented by the following diagram

\begin{figure}[h]
\begin{center}
\includegraphics[width=1.2in]{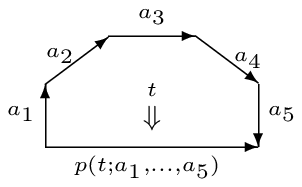}
\end{center}
\end{figure}

Given objects $x_i\in \ob(\mathcal{C})$, for $i=0,1,\dots,8$,
morphisms $a_i\in
\mo (x_{i-1},x_i)$, $i=1,2,\dots,8$, $t_1,t_3\in \mathcal{O}(3)$,
$t_2\in
\mathcal{O}(2)$ and $t\in
\mathcal{O}(4)$, then the axiom from Definition \ref{def5} is
represented by the following commutative diagram.

\begin{figure}[h]
\begin{center}
\includegraphics[width=4.5in]{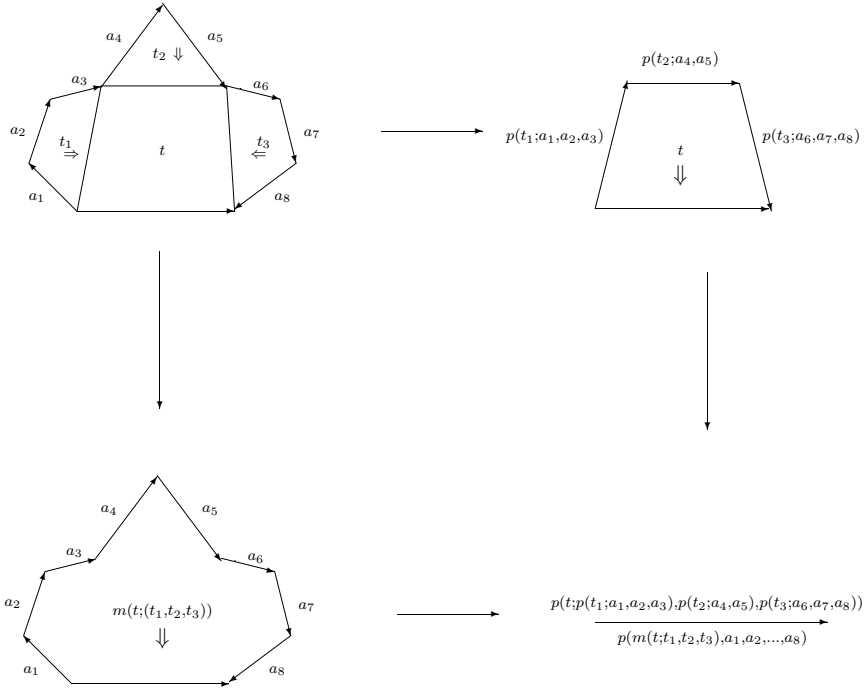}
\caption{Pictorial representation of the associative axiom of Definition
\ref{def5}. \label{fig:one}}
\end{center}
\end{figure}

Given a nonsymmetric operad $\mathcal{O}$, we define the category
$\ocat_k$ as follows:
\begin{itemize}

\item {$\ob(\ocat_k)= \mbox{small } \mathcal{O}$-categories.}

\item {Morphism $\mo_{\mathcal{O}}(\mathcal{C},\mathcal{D})$
from $\mathcal{O}$-category $\mathcal{C}$
to $\mathcal{O}$-category  $\mathcal{D}$ are functors from
$\mathcal{C}$ to $\mathcal{D}$ such that
$$F(p(t;a_1,\dots,a_n))=p(t;F(a_1),\dots,F(a_n)),$$ for given
objects $x_0,x_1,\dots,x_n$, morphisms  $a_i
\in\ho_{\mathcal{C}}(x_{i-1},x_1)$, $F(a_i)\in
\ho_{\mathcal{D}}(F(x_{i-1});F(x_i))$ and $t\in \mathcal{O}(n)$.
We call such a functor $F$ an
$\mathcal{O}$-functor.}
\end{itemize}

\section{Polya Functors}\label{fp}
Given an $\mathcal{O}$-category $\mathcal{C}$,
$\aut^{1}(\mathcal{C})\subset
\fu_\mathcal{O}(\mathcal{C},\mathcal{C})$ denotes the collection
of invertible $\mathcal{O}$-functors identical on objects. Let $G$
be a compact topological group. A $G$-action on an
$\mathcal{O}$-category $\mathcal{C}$ is a representation $\rho: G
\longrightarrow \aut^{1}(\mathcal{C}).$ It is defined by a
collection of actions $\rho_{x,y}: G \longrightarrow \mbox{{\rm
GL}}(\ho_\mathcal{C}(x,y))$ such that
$$\rho_{x_0,x_m}(g)(p(t;a_1,\dots,a_m))=p(t;\rho_{x_0,x_1}(g)(a_1),\dots, \rho_{x_{n-1},x_{m}}(g)(a_m))$$
for objects $x_0,x_1,\dots,x_n\in\ob(\mathcal{C})$,
$a_i\in\ho(x_{i-1},x_i)$, for all $i\in [m]$ and $g\in G$. Abusing
notation we shall write $ga$  instead of $\rho_{x,y}(g)(a)$ where
$x,y\in\ob(\mathcal{C})$ and $a\in\ho(x,y)$ . We define
$\ocat_k(G)$ to be the category of all linear
$\mathcal{O}$-categories provided with $G$ actions. Morphisms $F$
from $\mathcal{O}$-category $\mathcal{C}$ to
$\mathcal{O}$-category $\mathcal{D}$ are $G$-equivariant
$\mathcal{O}$-functors $F$ from $\mathcal{C}$ into $\mathcal{D}$,
i.e., $F(ga)=gF(a)$, for all $a\in \mo (x,y)$, $g\in G$,  where
$x,y\in \ob(\mathcal{C})$. We define
$$\begin{array}{cccc}
  \inv: & \ocat_k(G) & \longrightarrow & \ocat_k \\
  \mbox{ } & \mathcal{C} & \longmapsto & \mathcal{C}^{G}
\end{array} $$
as follows: $\ob({\mathcal{C}}^{G})=\ob(\mathcal{C})$ ,
$\mo_{\mathcal{C}^{G}}(x,y)=\mo_{{\mathcal{C}}}(x,y)^{G}$.
$\mathcal{C}^{G}$ is a $\mathcal{O}$-category since
$$g p(t;a_1,\dots,a_m)=p(t;g a_1,\dots,g a_m)=p(t;a_1,\dots,a_m).$$
We define
$$\begin{array}{cccc}
  \coinv: & \ocat_k(G) & \longrightarrow & \ocat_k \\
  \mbox{ } & \mathcal{C} & \longmapsto & {\mathcal{C}}_G
\end{array} $$
as follows:  $\ob({\mathcal{C}}_G)=\ob(\mathcal{C})$,
$\mo_{{\mathcal{C}}_G}(x,y)=\mo_\mathcal{C}(x,y)_G$ and

\begin{equation}\label{fone}{\displaystyle \overline{p}(t;\overline{a_1},\overline{a_2},
\dots ,\overline{a_m}) =\frac{1}{\mbox{vol}(G)^{m-1}}
\overline{\int_{G^{ (n-1)}} p(t;a_1,g_{2}a_{2},g_{3}a_{3}, \dots ,
g_{m}a_{m})\dif g_{2}\dif g_{3} \dots \dif g_{m}}}\end{equation}
\begin{thm}
There is a natural isomorphism $s:\coinv \longrightarrow \inv$
i.e., for $\mathcal{C}\in\ob(\ocat_k(G))$ we are given an
isomorphism $s_{\mathcal{C}}: \coinv(\mathcal{C}) \longrightarrow
\inv(\mathcal{C})$ such that $\inv(F)\circ s_{\mathcal{C}}=s_{\mathcal{D}}\circ
\coinv(F)$, for all functors $F:\mathcal{C}\longrightarrow
\mathcal{D}$ in the category $\ocat_k(G)$.
\end{thm}
\begin{proof} Given a category $\mathcal{C}$ and objects $x,y\in
\ob(\mathcal{C})$, let
$s_{x,y}:\ho_\mathcal{C}(x,y)\longrightarrow
\ho_\mathcal{C}(x,y)^{G}$  the symmetrization map defined in
the Section $\ref{sinco}$. By Proposition $\ref{inco}$, it induces
an isomorphism of vector spaces
$$s_{x,y}:\ho_\mathcal{C}(x,y)_G\longrightarrow \ho_\mathcal{C}(x,y)^{G}.$$
It remains to check that the following diagram is commutative

\[ \xymatrix{
&  \mathcal{O}(k)\otimes\ho_{\mathcal{C}}(x_0,x_1)_G\otimes \dots
\otimes \ho_{\mathcal{C}}(x_{m-1},x_m)_G \ar[dr]^{{\rm id}\otimes
s}\ar[ld]_{p(t,-)}& \\
\ho_\mathcal{C}(x_0,x_m)_{G} \ar[dr]_{s_{x_0,x_m}}& &
\hspace{-5.5cm}\mathcal{O}(k)\otimes\ho_{\mathcal{C}}(x_0,x_1)^{G}\otimes
\dots
\otimes \ho_{\mathcal{C}}(x_{m-1},x_m)^{G} \ar[ld]^{p(t,-)}\\
& \ho_\mathcal{C}(x_0,x_m)^{G} & }
\]
where $s= s_{x_0,x_1}\otimes \dots \otimes  s_{x_{m-1},x_m}$.
Clearly,
$$p(t;-) {\rm {id}}\otimes S=\frac{1}{\mbox{vol}(G)^{m}}\int_{G^{m}}
 p(t;g_1 a_1,\dots, g_m a_m) \dif g_{1}\dif g_{2} \dots \dif g_{m}.$$
On the other hand,
\begin{eqnarray*}
s_{x_0, x_m}p(t; \overline{a_1},\dots, \overline{a_m})&=&
\frac{1}{\mbox{vol}(G)}\int_{G}
g_1\overline{p(t; a_1,\dots, a_m)}\dif g_{1}\\
\mbox{} &\mbox{}& \mbox{}\\
\mbox{} &=&\frac{1}{\mbox{vol}(G)^{m}}\int_{G^{m}} g_1 p(t; a_1,
h_2a_2,\dots, h_m a_m) \dif
g_{1}\dif h_{2} \dots \dif h_{m}\\
\mbox{} &\mbox{}& \mbox{}\\
\mbox{} &=&\frac{1}{\mbox{vol}(G)^{m}}\int_{G^{m}}  p(t; g_1 a_1,
g_1h_2a_2,\dots, g_1h_m a_m) \dif
g_{1}\dif h_{2} \dots \dif h_{m}\\
\mbox{} &\mbox{}& \mbox{}\\
\mbox{} &=&\frac{1}{\mbox{vol}(G)^{m}}\int_{G^{m}}  p(t;g_1
a_1,\dots, g_m a_m) \dif g_{1}\dif g_{2} \dots \dif g_{m}
\end{eqnarray*}
making the change of variables $g_2=g_1 h_2, \dots, g_m=g_1 h_m$.
\end{proof}
Notice that if $\mathcal{O}$ is an operad then $\mathcal{O}\otimes
\mathcal{O}$ is naturally an operad with $(\mathcal{O}\otimes
\mathcal{O})(n)=\mathcal{O}(n)\otimes \mathcal{O}(n)$.
\begin{defi}
A Hopf operad is an operad together with an operad morphism
$\Delta:\mathcal{O}\longrightarrow \mathcal{O}\otimes
\mathcal{O}$.
\end{defi}
We have the following lemma
\begin{lema}
If $\mathcal{O}$ is a Hopf operad then the category of
$\mathcal{O}$-algebras is monoidal.
\end{lema}
\begin{proof}If $A$ and $B$ are $\mathcal{O}$-algebras then
$A\otimes B$ is also an $\mathcal{O}$-algebra, as the following
diagram shows
$$\mathcal{O}(n)\otimes (A\otimes B)^{\otimes n}\longrightarrow
(\mathcal{O}(n)\otimes A^{\otimes n}) \otimes
(\mathcal{O}(n)\otimes B^{\otimes n})\longrightarrow A\otimes B.$$
\end{proof}
We now construct a partial monoidal structure on $\ocat$.
\begin{defi}
Let $\mathcal{O}$ be a Hopf operad, given $\mathcal{O}$-categories
$\mathcal{C}$ and $\mathcal{D}$ such that
$\ob(\mathcal{C})=\ob(\mathcal{D})=X$, we define the tensor
product category $\mathcal{C}\otimes\mathcal{D}$ as follows
\begin{itemize}
\item{$\ob(\mathcal{C}\otimes \mathcal{D})=X$.}
\item{$\mo_{{\mathcal{C}}\otimes {\mathcal{D}}}
(x,y)=\mo_{\mathcal{C}}(x,y)\otimes \mo_{\mathcal{D}}(x,y)$.}
\item{$p(t; a_1\otimes b_1,\dots,a_n\otimes b_n)=\sum
p(t_{(1)};a_1,\dots,a_n)\otimes p(t_{(2)};b_1,\dots,b_n)$,  where
$\Delta(t)=\sum t_{(1)}\otimes t_{(2)}$ using Swedler notation.}
\end{itemize}
\end{defi}
Recall the well-known definition. Given a pair of groups $G$ and
$K\subset S_n$ the semidirect product $G^{n}\rtimes K$ is the set
$G^{n}\rtimes K=\{(g,a): g\in G^{n}, a\in K\}$ provided with the
product $(g,a)(h,b)=(ga(h),ab)$ where $g,h\in G^{n}$, $a,b\in K$
and if $h=(h_1,\dots,h_n)$ then \
$a(h_1,\dots,h_n)=(h_{a^{-1}(1)},\dots,h_{a^{-1}(n)}).$ The
following result is obvious
\begin{lema}
\begin{enumerate}
\item{$\ob(\mathcal{C}^{\otimes n})=\ob(\mathcal{C})$, $\mo_{\mathcal{C}^{\otimes n}}(x,y)=
\mo_{\mathcal{C}}(x,y)^{\otimes n}$.}
\item{If $G$ acts on $\mathcal{C}$ and $K\subset S_n$ then $G^{n}\rtimes K$
acts on $\mathcal{C}^{\otimes n}$.}
\end{enumerate}
\end{lema}
Assume that $\mathcal{C}$ is a category over a non-symmetric Hopf
operad $\mathcal{O}$. Let $G$ be a compact topological group, and
$K$ a subgroup of $S_n$. We construct functor $P_{G,K}$ which we
call the {\it Polya functor} of type $G,K$
$$P_{G,K}:
\ocat_k(G) \longrightarrow \ocat_k$$ as follows:

\begin{itemize}
\item{On objects: Given $\mathcal{C}\in \ob(\ocat_k(G))$, then
$P_{G,K}(\mathcal{C})\in \ob(\ocat_k)$ is the category
$$P_{G,K}(\mathcal{C}) = \mathcal{C}^{\otimes{n}}/G^{n}\rtimes K.$$ Explicitly:}
\item{ ${\ob}(P_{G,K}(\mathcal{C}))=\ob(\mathcal{C}^{\otimes
n})=\ob(\mathcal{C})$, and for given objects $x,y\in
{\ob}(P_{G,K}(\mathcal{C})),$
$${\mo}_{P_{G,K}(\mathcal{C})}(x,y)=(\mo_{\mathcal{C}}(x,y)^{\otimes{n}})/{G^{n}\rtimes K}.$$}
\item{Identity: $\id_x \in
\mo_{P_{G,K}(\mathcal{C})}(x,x)=\overline{{\id}_x^{\otimes n} }\in
({\mo_{\mathcal{C}}}(x,x)^{\otimes{n}})/{G^{n}\rtimes K}$.}
\item{Composition: Given  $x_0,\dots,x_m\in\ob(P_{G,K}(\mathcal{C}))$ and morphisms $\overline{a_i}\in
{\mo}_{P_{G,K}(\mathcal{C})}(x_{i-1},x_i)$, for $i=1,\dots,m$,  we
have the following
$$p(t;\overline{a_1},\dots, \overline{a_m})=\frac{1}{((\sharp(G)^{n}) \sharp(K))^{m-1}}\sum_{(g,s)\in (G^{n}\rtimes
K)^{m-1}}\overline{p(t;a_1,(g_2,s_2)a_2,\dots,(g_m,s_m)a_m)}\ \
\mbox{if $G$ is finite}, $$
$$ p(t;\overline{a_1},\dots, \overline{a_m})=\frac{1}{(\vol(G^{n}) \sharp(K))^{m-1}}{\des\sum_{s\in K^{m-1}}
 \int_{g\in (G^{n})^{m-1}}
\!\!\overline{p(t;a_1,(g_2,s_2)a_2,\dots,(g_m,s_m)a_m)} \dif g_2\dots \dif
g_m,} $$ if $G$ is compact and $g_i\in G^{n}$.}

\item{On morphisms: each  functor $\xymatrix @1{\mathcal{C}
\ar[r]^{\alpha} & \mathcal{D}}$ induces a functor $\xymatrix
@1{\mathcal{C}^{\otimes n} \ar[r]^{\alpha^{\otimes n} } &
\mathcal{D}^{\otimes n}}$. This functor descends to a well defined
functor
\[
\xymatrix @1{P_{G,K}(\mathcal{C})=\mathcal{C}^{\otimes{n}}/G^{n}\rtimes K
\ar[r] & \mathcal{D}^{\otimes{n}}/G^{n}\rtimes
K=P_{G,K}(\mathcal{D})}.
\]}
\end{itemize}

\begin{exa}
Consider the non-symmetric operad $\mathcal{ASS}$ given by
$\mathcal{ASS}(n)=k$, for all $n\geq 0$.
$\mathcal{ASS}$-categories are categories in the usual sense. This
example will be applied in Section \ref{sct} to introduce explicit
formulae for the composition of morphisms in the Schur categories
of various types.
\end{exa}

\begin{defi}
Let $\mathcal{O}$ be an operad. An $\mathcal{O}$-algebra is a pair
$(A,\gamma)$ where $A$ is a vector space and\\
$\gamma:\mathcal{O}\longrightarrow \en_A$ is an operad morphism.
\end{defi}

Notice that an $\mathcal{O}$-algebra may be regarded as an
$\mathcal{O}$-category $\mathcal{C}$ with one object by setting
$A=\mo_{\mathcal{C}}(1,1)$. We denote by $ \oalg_k$ the category
of $\mathcal{O}$-algebras and by $\oalg_k(G)$ the category of
$\mathcal{O}$-algebras provided with a $G$ action. We have two
naturally isomorphic functors
$$\begin{array}{cccc}
  \inv: & \oalg_k(G) & \longrightarrow & \oalg_k \\
  \mbox{ } & A & \longmapsto & A^{G}
\end{array} $$
$$\begin{array}{cccc}
  \coinv: & \oalg_k(G) & \longrightarrow & \oalg_k \\
  \mbox{ } & A & \longmapsto & A_G
\end{array} $$
Let $\mathcal{O}$ be a Hopf operad. Assume that $A$ is an
$\mathcal{O}$-algebra. Let $G$ be a compact
 topological group and $K\subset S_n$. We have a functor
$$\begin{array}{cccc}
  P_{G,K}:
 & \oalg_k(G) & \longrightarrow &\oalg_k  \\
  & A & \longmapsto & A^{\otimes n}/G^{n}\rtimes K \\
\end{array}$$

The $\mathcal{O}$-algebra structure  on $A^{\otimes
n}/G^{n}\rtimes K $ is given for any $a_1,\dots,a_m\in A$,
$t\in\mathcal{O}(m)$ by

\begin{equation}\label{fdp}p(t;\overline{a_1},\dots, \overline{a_m})=\frac{1}{((\sharp(G)^{n}) \sharp(K))^{m-1}}\sum_{(g,s)\in (G^{n}\rtimes
K)^{m-1}}\overline{p(t;a_1,(g_2,s_2)a_2,\dots,(g_m,s_m)a_m)}
\end{equation}
 if
$G$ is finite. If $G$ is compact taking $g_i\in G^{n}$ we have
that
$$ p(t;\overline{a_1},\dots, \overline{a_m})=\frac{1}{(\vol(G^{n}) \sharp(K))^{m-1}}{\des\sum_{s\in K^{m-1}}
 \int_{g\in (G^{n})^{m-1}}
\!\!\overline{p(t;a_1,(g_2,s_2)a_2,\dots,(g_m,s_m)a_m)} \dif g_2\dots \dif
g_m}. $$
\begin{cor}\label{prod2} Let $A$ be an $\mathcal{O}$-algebra,  $G$ a finite group acting by algebra automorphism on $A$. Take
$K=\{\id\}$. The following identity hold in $P_G(A)$
\begin{equation}
\overline{a}\ \overline{b}=\frac{1}{\sharp(G)}\sum_{g\in G}
\overline{a (gb)}
\end{equation}
where $a,b\in A$.
\end{cor}
We remark that Joyal theory of analytic functors, see
$\cite{Joy1}$, $\cite{Joy2}$, may be extended from the context of
$k$-vector spaces to $\mathcal{O}$-algebras by defining a functor
$F$ which sends a family $A=\{A_n\}_{n\geq 0}$ of
$\mathcal{O}$-algebras provided with right $S_n$ actions, into the
functor $F_A$ from the $\oalg$ into $\oalg$ given as follows
$$\begin{array}{ccccc}
  F: & \fu(S^{op},\oalg ) & \longrightarrow & \fu(\oalg,\oalg)& \mbox{ } \\
   \mbox{ } & A=\{A_n\}_{n{\geq 0}} & \longmapsto  & F_{A}(B)={\des \bigoplus_{n\geq 0 }
   (A_n\otimes B^{\otimes n})_{S_n}}& \mbox{for all}\ \ B\in\ob(\oalg)\\
\end{array}$$

\begin{exa}
Consider the operad $Ass$ given by $Ass(n)=k[S_n]$, for $n\geq 0$.
$Ass$-algebras are the same as associative algebras. If we take
the family $k=\{k\}_n$, $n\geq 0$, provided with the trivial $S_n$
action, then $F_k(A)=\sym(A)$
\end{exa}

The following remarks justify our choice of name for the Polya
functors. Let $K\subset S_n$ be a permutation group. For $k\in K$,
let $b_s(k)$ be the number of cycles of $k$ of length $s$. The
{\em cycle index polynomial} of $K\subset S_n$, is the polynomial
in $n$ variables $x_1,\dots,x_n$
$$P_K(x_1,\dots,x_n)=\frac{1}{\sharp(K)}\sum_{k\in K} x_1^{b_1(k)} x_2^{b_2(k)}\dots x_n^{b_{n}(k)}.$$
Let $A$ be a finite dimensional $\mathcal{O}$-algebra and $X$ a
basis for $A$. Using Polya theory $\cite{KW}$, we can compute the
dimension of the $\mathcal{O}$-algebra $(A^{\otimes n})_K$ as
follows:
\begin{eqnarray*}
\dim ((A^{\otimes n})_K)&=&\sharp(X^{n}/K)\\
\mbox{}&=&P_K(\sharp(X),\sharp(X),\dots,\sharp(X))\\
\mbox{}&=&P_K(\dim A,\dim A,\dots,\dim A).
\end{eqnarray*}

\subsection{General multiplication rule }

For each $K\subset S_n$ consider the Polya functor $P_K:k\mbox{-}
\alg\longrightarrow k\mbox{-}\alg$ from the category of
associative $k$-algebras into itself defined on objects as
follows: if $A$ is a $k$-algebra, then $P_K(A)$  denotes the
algebra whose underlying vector space is $$P_K(A)= (A^{\otimes
n})/ \langle a_1\otimes \dots \otimes a_n -a_{\sigma^{-1}
(1)}\otimes \dots \otimes a_{\sigma^{-1}(n)}: a_i\in A, \sigma\in
K\rangle .$$ Our next theorem provides the rule for the product of
$m$ elements in $P_K(A)$.
\begin{thm}\label{pfu} For any $a_{ij}\in A$ the following
identities holds in $P_K(A)$
\begin{equation}\label{PF}
\sharp(K)^{m-1}\prod_{i=1}^{m}\left(\overline{\bigotimes_{j=1}^{n}
a_{ij}}\right)  = \sum_{\sigma \in \{\id\} \times
K^{m-1}}\overline{\bigotimes_{j=1}^{n}\left( \prod_{i=1}^{m} a_{i
\sigma^{-1}_i(j)}\right)}  \end{equation}
\end{thm}
\begin{proof}
It is follows from formula (\ref{fdp}), taking $Ass$ as the
underlying operad and setting $G=\{\id\}$.
\end{proof}

Theorem $\ref{pfu}$ implies the following

\begin{prop}\label{apf}
Let $A$ be an algebra provided with a basis  $\{e_s| s\in[r]\}$.
Assume that $e_s e_t=c(k,s,t) e_k$ (sum over $k$), for all
$s,t\in[r]$. For any given  $a=(a_{ij})\in \mbox{M}_{n\times
m}([r])$ the following identity holds in $P_K(A)$
$$(n!)^{m-1}{\des \prod_{i=1}^{m}\left(\overline{\bigotimes_{j=1}^{n}
e_{a_{ij}}}\right)}=\sum_{\sigma, \alpha}
\left(\prod_{j=1}^{n}c(\alpha^{j},a,\sigma)\right)
\overline{\bigotimes_{j=1}^{n} e_{\alpha_{ m-1}^{j}}}, $$ where
the sum runs over all $\sigma\in \{\id\} \times K^{m-1}$,
$\alpha=(\alpha_i^{j})\in{\mbox{M}}_{(m-1)\times n}([r])$, and\\
$ c(\alpha^{j},a,\sigma)= {\des c(\alpha_1^{j}, a_{1
\sigma^{-1}_1(j)}, a_{2 \sigma^{-1}_2(j)}) c(\alpha_2^{j},
\alpha_1^{j}, a_{3 \sigma^{-1}_3(j)}) \dots c(\alpha_{m-1}^{j},
\alpha_{m-2}^{j}, a_{m\sigma^{-1}_m(j)})}$.
\end{prop}

\begin{proof}
\begin{eqnarray*}
(n!)^{m-1}{\des
\prod_{i=1}^{m}\left(\overline{\bigotimes_{j=1}^{n}
e_{a_{ij}}}\right)}&=& {\des \sum_{\sigma\in \{ \id \}\times
G^{m-1}} \overline{\bigotimes_{j=1}^{n}\left( \prod_{i=1}^{m}
e_{a_{i \sigma_i^{-1}(j)}}\right)}}\\
\mbox{} &\mbox{}&\mbox{}\\
\mbox{}&=& {\des \sum_{\alpha, \sigma}\overline{
\bigotimes_{j=1}^{n}
c(\alpha,a,\sigma) e_{\alpha_{m-1}}}}\\
\mbox{} &\mbox{}&\mbox{}\\
\mbox{}&=& {\des \sum_{\sigma, \alpha}
\left(\prod_{j=1}^{n}c(\alpha^{j},a,\sigma)\right)
\overline{\bigotimes_{j=1}^{n} e_{\alpha_{ m-1}^{j}}}}
\end{eqnarray*}

\end{proof}
Let us consider a non-symmetric Hopf operad $\mathcal{O}$. Assume
that a basis $p_m^{t}$, $t\in[k_m]$, for $\mathcal{O}(m)$ is given
for each $m\in\mathbb{N}$. Moreover let us assume that
$\Delta(p_m^{t})=p_m^{t}\otimes\dots\otimes p_m^{t}$, then the
next proposition follows from formula $(\ref{PF})$.

\begin{prop}
Let $A$ be an $\mathcal{O}$-algebra provided with a basis  $\{e_s|
s\in[r]\}$, and let $\mathcal{O}(n)=\langle p_m^{t}:
t\in[r_m]\rangle$. Assume that
$p_m^{t}(e_{s_1},\dots,e_{s_m})=c(t,k,m,s_1,\dots,s_n) e_k$, (sum
over $k$). For any given $a=(a_{ij})\in\mbox{M}_{n\times m}([r])$,
the following identity holds in $P_K(A)$
$$p_m^{t}( \overline{\otimes_{j=1}^{n}
e_{a_{1j}}},\dots,\overline{\otimes_{j=1}^{n} e_{a_{mj}}})=
\sum_{\sigma,u}\prod_{j=1}^{n}c(t,u_j,m,a_{1
\sigma^{-1}_1(j)},\dots,a_{m \sigma^{-1}_m(j)}) \overline{
\bigotimes_{j=1}^{n} e_{u_j}}$$ where the sum runs over all
$\sigma\in \{\id\} \times K^{m-1}$ and $u\in[r]^{n}$.
\end{prop}

\section{Symmetric power of a supercategory.}\label{sct}

Let us consider Polya functor $P_{G,K}:\cat_k(G)\longrightarrow
\cat_k$ for the case $G=\id$, $K=S_n$, i.e, we consider for each
$n\in\mathbb{N}$ the functor
$$\sym^{n}: \cat_k\longrightarrow \cat_k$$
Recall that a supercategory is a category over the category {\bf
{Supervect}} of $\mathbb{Z}_2$-graded vector spaces with the
Koszul rule of signs. Functor $\sym^{n}$ may be applied to
supercategories as well. The next result provides formula for the
composition of morphisms in the symmetric powers of a
supercategory. We use the notation $a_1\dots
a_n=\overline{a_1\otimes \dots \otimes a_n}\in
\sym^{n}(\ho_{\mathcal{C}}(x,y))$, for morphisms $a_1,\dots,a_n\in
\ho_{\mathcal{C}}(x,y)$.

\begin{prop}\label{psym}
Let $\mathcal{C}$ be a supercategory and let $a_1,\dots,a_n \in
\mo_{\mathcal{C}}(x,y)$ and $b_1,\dots,b_n \in
\mo_{\mathcal{C}}(y,z)$. In the supercategory
$\sym^{n}(\mathcal{C})$ the compositions of morphisms is given by
$$(a_{1} a_{2}  \dots a_{n})(  b_{1} b_{2}
\dots b_{n}) = \frac{1}{n!}\sum_{\sigma \in S_n}\sgn(a,b,\sigma)
(a_{1}b_{\sigma^{-1}(1)})( a_{2}b_{\sigma^{-1}(2)})
  \dots (a_{n}b_{\sigma^{-1}(n)})$$
where $\sgn(a,b,\sigma)=(-1)^{e}$ and $e=e(a,b,\sigma)=
\sum_{i>j}\overline{a_{i}}\overline{b_{\sigma^{-1}(j)}} +
\sum_{\sigma(i)>\sigma(j)}\overline{b_i}
\overline{b_j}.$
\end{prop}

\subsection{Schur categories}

Let $k$ be a field of characteristic $0$,
$m=(m_1,\dots,m_k)\in\mathbb{N}^{k}$,
$\mathbb{Z}_m=\mathbb{Z}_{m_1}\times \dots \times
\mathbb{Z}_{m_k}$ and $n\in\mathbb{N}$. We define the Schur
supercategory of type $(m,n)$  as follows:
$$\ob(S(m,n))=\mbox{finite dimesional $k$-supervector spaces}.$$
$$\mo_{S(m,n)}(V,W)=(\ho_k(V^{\oplus \mathbb{Z}_m}, W^{\oplus \mathbb{Z}_m})^{\otimes n})_{\mathbb{Z}_m^{n}\rtimes S_n}.$$
$\mathbb{Z}_m^{n}\rtimes S_n$ acts on $(\ho_k(V^{\oplus
\mathbb{Z}_m}, W^{\oplus \mathbb{Z}_m})^{\otimes n})$ as follows
$$\begin{array}{ccc}
  \mathbb{Z}_m^{n}\rtimes S_n \times (\ho_k(V^{\oplus \mathbb{Z}_m}, W^{\oplus \mathbb{Z}_m})^{\otimes n})
  & \longrightarrow & (\ho_k(V^{\oplus \mathbb{Z}_m}, W^{\oplus \mathbb{Z}_m})^{\otimes n}) \\
   &  &  \\
  ((c_1,\dots,c_n),\sigma)(E_{r_1s_1}^{t_1u_1}\dots E_{r_n s_n}^{t_n u_n})
  &  \longmapsto &{\des (E_{r_{\sigma(1)}(s_1+c_1)}^{t_{\sigma(1)}(u_1+c_1)}
  \dots E_{r_{\sigma(n)}(s_n +c_n)}^{t_{\sigma
  (n)}(u_n+c_n)})}
\end{array}$$
where $E(V,W)_{ij}^{kl}$ are the elementary linear transformation
in $(\ho_k(V^{\oplus
\mathbb{Z}_m}, W^{\oplus
\mathbb{Z}_m})$, $i\in[\dim V]$,\\ $k\in[\dim
W]$ and $j,l\in \mathbb{Z}_m^{k}$. We apply Polya functor to
obtain explicit formula for the composition  rule
$$\mo(V,W)\otimes \mo(W,Z) \longrightarrow \mo(V,Z).$$

\begin{thm}\label{sc} For any given $M=m_1\dots m_k$,
$i,t\in [\dim W]$, $k\in [\dim Z]$, $r\in[\dim V]$, and
$j,l,s,u\in\mathbb{Z}_m^{n}$, we have
$$(E(V,W)_{i_1 j_1}^{k_1l_1}\dots E(V,W)_{i_nj_n}^{k_nl_n})
(E(W,Z)_{r_1s_1}^{t_1 u_1}\dots E(W,Z)_{r_ns_n}^{t_nu_n})=$$
$$\frac{1}{M^{n} n!}\sum_{\begin{array}{c}
  \sigma \in S_n \\
  t_{\sigma(a)}=i_a \\
\end{array}} \sgn(\sigma,i,k,r,t) E(V,Z)_{r_{\sigma(1)}(s_1+j_1-u_1)}^{k_1 l_1}\dots
E(V,Z)_{r_{\sigma(n)}(s_n+j_n-u_n)}^{k_n l_n}$$ where
$\sgn(\sigma,i,k,r,t)=(-1)^{e}$ and
$$e={\des \sum_{i>j}
(i_i+k_i)(r_{\sigma^{-1}(i)}+ t_{\sigma^{-1}(i)})}+{\des
\sum_{\sigma(i)>\sigma(j)}(r_{\sigma^{-1}(i)}+
t_{\sigma^{-1}(i)})(r_{\sigma^{-1}(j)}+ t_{\sigma^{-1}(j)})}.$$
\end{thm}
\begin{proof} Straightforward using Polya functor and Proposition
\ref{psym}.\end{proof}

\begin{figure}[ht]
\begin{center}
\epsfxsize 10cm \epsfbox{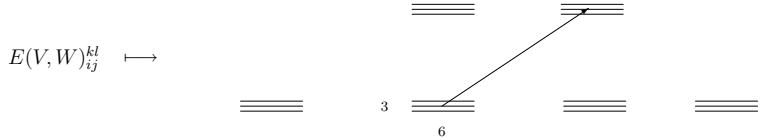} \caption{Representation of
elementary transformation. \label{fig:grafrep}}
\end{center}
\end{figure}

We now develop a graphical notation that make transparent the
meaning of Theorem $\ref{sc}$. Let us assume that $n=4$, $k=4$,
$m_1=6$, $m_2=3$, $\dim (V)=4$, $\dim (W)=2$ and $\dim (Z)=3$. We
represent an elementary linear transformation $E(V,W)$ as in
Figure $\ref{fig:grafrep}$

Notice that each block corresponds to $\mathbb{Z}_6\times
\mathbb{Z}_3$, $\mathbb{Z}_6$ acting horizontally, and
$\mathbb{Z}_3$ acting vertically. The number of blocks in the
bottom row is $\dim(V)$, and the number of blocks in the top row
is $\dim(W)$. Elements of $(\ho(V^{\oplus \mathbb{Z}_m},W^{\oplus
\mathbb{Z}_m})^{\otimes n})_{(\mathbb{Z}_{6}\times
\mathbb{Z}_{3})^{4}{\rtimes S_4}}$ are depicted by four
non-numbered arrows, and similarly  for elements of
$(\ho(W^{\oplus \mathbb{Z}_m},Z^{\oplus \mathbb{Z}_m})^{\otimes
n})_{(\mathbb{Z}_{6}\times \mathbb{Z}_{3})^{4}{\rtimes S_4}}$.
Composition is obtained as follows

\begin{figure}[ht]
\begin{center}
\includegraphics[width=5in]{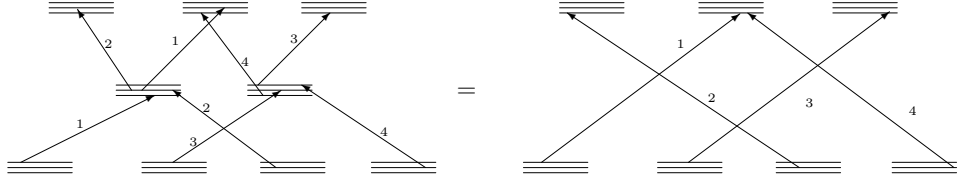}
\caption{Example of composition. \label{fig:two}}
\end{center}
\end{figure}

\begin{itemize}
\item{Fix an arbitrary enumeration of the arrows in
$(\ho(V^{\oplus \mathbb{Z}_m},W^{\oplus \mathbb{Z}_m})^{\otimes
n})_{(\mathbb{Z}_{6}\times \mathbb{Z}_{3})^{4}{\rtimes S_4}}$.}
\item {Sum over all possible enumerations of the arrows in
$(\ho(W^{\oplus \mathbb{Z}_m},Z^{\oplus \mathbb{Z}_m})^{\otimes
n})_{(\mathbb{Z}_{6}\times \mathbb{Z}_{3})^{4}{\rtimes S_4}}$.}
\item {Stacks arrows from $(\ho(V^{\oplus \mathbb{Z}_m},W^{\oplus
\mathbb{Z}_m})^{\otimes n})_{(\mathbb{Z}_{6}\times
\mathbb{Z}_{3})^{4}{\rtimes S_4}}$ to arrows on  $(\ho(W^{\oplus
\mathbb{Z}_m},Z^{\oplus \mathbb{Z}_m})^{\otimes
n})_{(\mathbb{Z}_{6}\times \mathbb{Z}_{3})^{4}{\rtimes S_4}}$
taking care of enumeration and using the $\mathbb{Z}_2\times
\mathbb{Z}_6$ symmetry. }
\end{itemize}
Notice that composition is interesting in that no-touching arrows
may nevertheless be composed (due to the $\mathbb{Z}_{6}\times
\mathbb{Z}_{3}$ symmetry) as shown in Figure $\ref{fig:two}$.

\begin{defi} The Schur superalgebra
of type $(\sdim V,m)$ is given by $\ho_{S_n}(V^{\otimes
n},V^{\otimes n})$, where $\sdim$ denotes the superdimension of a
supervector space.
\end{defi}
See $\cite{JG}$ for more on Schur algebras.
\begin{cor}
For $m=1$ and $V=W$, $\mo_{S(m,n)}(V,V)$ is the Schur
superalgebra\ $\schur(\sdim V,m)$ of type $(\sdim V,m)$.
\end{cor}
\begin{proof}   $(\ho(V,V)^{\otimes n})_{S_n}\cong (\ho(V,V)^{\otimes
n})^{S_n}\cong \ho_{S_n}(V^{\otimes n},V^{\otimes n})$.
\end{proof}

\section{Classical symmetric functions}

In this section we study classical symmetric functions by means of
the Polya functor. We provide a fairly elementary interpretation
of the symmetric functions in terms of the symmetric powers of the
monoidal algebra associated to the additive monoid
$\mathbb{N}^{m}$. We also consider symmetric odd-functions as well
as symmetric Boolean algebras. Symmetric functions have been
studied from many points of view, see for example $\cite{Macd}$,
$\cite{GR}$, $\cite{VC}$ .

\subsection{Symmetric functions of Weyl type}
The classical Weyl groups of type $A_n, B_n,$ and $D_n$, are
$S_n$, $\mathbb{Z}_2^{n} \rtimes S_n$ and $\mathbb{Z}_2^{n-1}
\rtimes S_n$ respectively. These groups act on
$(\mathbb{R}^{m})^{n}$ as follows:

$$\begin{array}{ccc}
  S_n \times (\mathbb{R}^{m})^{n} & \longrightarrow & (\mathbb{R}^{m})^{n} \\
  (\sigma, (x_1,\dots, x_n)) & \longmapsto & (x_{\sigma^{-1}(1)},\dots, x_{\sigma^{-1}(n)}) \\
\end{array}$$

$$\begin{array}{ccc}
 (\mathbb{Z}_2^{n} \rtimes S_n) \times (\mathbb{R}^{m})^{n} & \longrightarrow & (\mathbb{R}^{m})^{n} \\
  ((t_1,\dots,t_n),\sigma )(x_1,\dots, x_n) & \longmapsto & (t_1x_{\sigma^{-1}(1)},\dots, t_n x_{\sigma^{-1}(n)}) \\
\end{array}$$
The group $D_n$ is regarded as a subgroup of $B_n$ as follows
$$\mathbb{Z}_2^{n-1}
\rtimes S_n=\{((t_1,\dots,t_n), \sigma)\in
\mathbb{Z}_2^{n} \rtimes S_n :  t_1 t_2\dots t_n=1
\}.$$

\begin{defi}
Fix $m\in \mathbb{N}$. The algebra of symmetric functions of type
$A_n$, $B_n$ and $D_n$ are given by
\begin{itemize}
\item{$\sym_{A_n}(m)=(\mathbb{C}[x_1,\dots,x_n])_{S_n}\cong
(\mathbb{C}[x_1,\dots,x_n])^{{S_n}}$, }
\item{$\sym_{B_n}(m)=
(\mathbb{C}[x_1,\dots,x_n])_{{\mathbb{Z}_2^{n}\rtimes S_n}}\cong
(\mathbb{C}[x_1,\dots,x_n])^{{{\mathbb{Z}_2^{n}\rtimes S_n}}}$,}
\item{ $\sym_{D_n}(m)=
(\mathbb{C}[x_1,\dots,x_n])_{{\mathbb{Z}_2^{n-1}\rtimes S_n}}\cong
(\mathbb{C}[x_1,\dots,x_n])^{{{\mathbb{Z}_2^{n-1}\rtimes S_n}}}$,}
\end{itemize}
where  $x_i=(x_{i1},\dots, x_{im})$, for $i=1,\dots,n$.
\end{defi}
The map  $\mathbb{C}[\mathbb{N}^{m}]^{\otimes n}  \longrightarrow
\mathbb{C}[x_1,\dots,x_n]$ given by $a_1\otimes\dots \otimes a_n  \longmapsto  x_1^{a_1}\dots x_n^{a_n}$
defines an isomorphism of algebras, where $a_i\in \mathbb{N}^{m}$
and $x_i^{a_i}=x_{i1}^{a_{i1}}\dots x_{im}^{a_{im}}$. We set
$\mathbb{N}_{e}^{m}=\{a\in\mathbb{N}^{m}: |a|\ \
\mbox{is even} \}$ and $\mathbb{N}_{o}^{m}=\{a\in\mathbb{N}^{m}:
|a|\ \ \mbox{is odd} \}$. We denote $X^{A}=x_1^{a_1}\dots
x_n^{a_n}$, for $A=(a_1,\dots,a_n)\in (\mathbb{N}^{m})^{n}$.

\begin{thm}[Classical symmetric functions]\label{csf}\quad

\begin{enumerate}
\item{The following is a commutative diagram
$$\begin{array}{ccccc}
  (\mathbb{C}[\mathbb{N}^{m}]^{\otimes n})_{S_n} & \supset & (\mathbb{C}[\mathbb{N}^{m}_e +
\mathbb{N}^{m}_o]^{\otimes n})_{S_n} & \supset & (\mathbb{C}[\mathbb{N}^{m}_e]^{\otimes n})_{S_n} \\
  \downarrow &  & \downarrow &  & \downarrow \\
  \sym_{A_n} (m) & \twoheadrightarrow &  \sym_{D_n}(m) & \twoheadrightarrow & \sym_{B_n}(m) \\
\end{array}$$ where the vertical arrows are isomorphisms.}
\item{ For $A,B\in (\mathbb{N}^{m})^{n}$, the product rule in $\sym_{A_n}(m)$ is given by
$$\overline{X^{A}}\ \  \overline{X^{B}}={\displaystyle \frac{1}{n!}\sum_{\sigma \in S_n}\overline{X^{A+\sigma(B)}}}$$
 and on $\sym_{B_n}(m)$ and\ \  $\sym_{D_n}(m)$
 by restriction.}
\end{enumerate}
\end{thm}
\begin{proof} The first row of the diagram above follows from the
isomorphism above after taken care of the $\mathbb{Z}_2^{n}$
(resp. $\mathbb{Z}_2^{n-1}$) symmetries for the groups $B_n$ and
$D_n$ respectively. It is clear that part $b$ implies the rest of
$a.$ We prove $b.$ Given $A,B\in(\mathbb{N}^{m})^{n}$, using
Proposition \ref{psym} for the product in $\sym_{A_n}(m)$, we
obtain
\begin{eqnarray*}
\overline{X^{A}}\ \  \overline{X^{B}} &=&\frac{1}{n!} \sum_{\sigma \in S_n} \overline{(x_1^{a_1}\dots
x_n^{a_n})
(x_{1}^{b_{\sigma^{-1}(1)}}\dots x_{n}^{b_{\sigma^{-1}(n)}})}\\
\mbox{} &\mbox{}& \mbox{}\\
\mbox{} &=&\frac{1}{n!}\sum_{\sigma \in S_n}\overline{X^{A+\sigma(B)}}.
\end{eqnarray*}
Now consider $A,B\in(\mathbb{N}_e^{m})^{n}$, the product in
$\sym_{B_n}(m)$ is given by
\begin{eqnarray*}
 \overline{X^{A}}\ \  \overline{X^{B}}&=&\frac{1}{2^{n}n!}\sum_{(t,\sigma) \in
\mathbb{Z}_2^{n}\rtimes S_n} \overline{ x_1^{a_1}\dots x_n^{a_n}
t_1^{b_1}x_{1}^{b_{\sigma^{-1}(1)}}\dots  t_n^{b_n}x_{n}^{b_{\sigma^{-1}(n)}}}\\
\mbox{} &\mbox{}& \mbox{}\\
\mbox{}&=&\frac{1}{2^{n}n!} \sum_{\sigma \in S_n}\left( \sum_{t \in
\mathbb{Z}_2^{n}}(-1)^{\displaystyle{\sum_{t_i=-1}
b_{i}}}\right)\overline{X^{A+\sigma(B)}}\\
\mbox{} &\mbox{}& \mbox{}\\
\mbox{} &=&\frac{1}{n!}\sum_{\pi \in S_n}\overline{X^{A+\sigma(B)}},\ \
\mbox{since}\ \ \sum_{t \in
\mathbb{Z}_2^{n}}(-1)^{\displaystyle{\sum_{t_i=-1}
b_{i}}}=2^{n}.
\end{eqnarray*}
Consider $A,B\in(\mathbb{N}_e^{n}+\mathbb{N}_o^{m})^{n}$, we
obtain
\begin{eqnarray*}
\overline{X^{A}}\ \  \overline{X^{B}} &=&\frac{1}{2^{n-1}n!}\sum_{(t,\sigma)
\in
\mathbb{Z}_2^{n}\rtimes S_n}  \overline{x_1^{a_1}\dots x_n^{a_n}
t_1^{b_1}x_{1}^{b_{\sigma^{-1}(1)}}\dots  t_n^{b_n}x_{n}^{b_{\sigma^{-1}(n)}}}\\
\mbox{} &\mbox{}& \mbox{}\\
\mbox{}&=&\frac{1}{2^{n-1}n!} \sum_{\pi \in S_n}\left( \sum_{\begin{array}{c}
  {\scriptstyle{t\in \mathbb{Z}_2^{n}}} \\
  {\scriptstyle{\prod t_i=1}} \\
\end{array}}
(-1)^{{\displaystyle \sum_{t_i=-1} b_{i}+b_n}}\right)\overline{X^{A+\sigma(B)}}\\
\mbox{} &\mbox{}& \mbox{}\\
\mbox{} &=&\frac{1}{n!}\sum_{\sigma \in S_n}\overline{X^{A+\sigma(B)}},\ \ \mbox{since}\ \ \sum_{\begin{array}{c}
  {\scriptstyle{t\in \mathbb{Z}_2^{n}}} \\
  {\scriptstyle{\prod t_i=1}} \\
\end{array}}(-1)^{{\displaystyle \sum_{t_i=-1} b_{i}+b_n}}=2^{n-1}
\end{eqnarray*}
\end{proof}

\subsection{Symmetric odd-functions}
Consider the alternating algebra $\bigwedge
[\theta_1,\dots,\theta_m]$, which we regard as the algebra of
functions on the purely odd super-space $\mathbb{R}^{0|m}$. A
basis for $\bigwedge [\theta_1,\dots,\theta_m]$ is given by
$\{\theta_{I}=\theta_{i_1}\dots \theta_{i_k}| \  I\subset [m] \}$.
The structural coefficients are given by $\theta_I
\theta_J=c(I,J)\theta_{I\cup J}$, where $c(I,J)= (-1)^{|\{ i\in I,
j\in J,\ i>j\}|}$, if $I\cap J=\emptyset$, and $0$ otherwise. The
algebra $\bigwedge [\theta_1,\dots,\theta_m]$ is
$\mathbb{Z}_2$-graded, with grading
$\overline{\theta_I}=\overline{I}=0$ if $|I|$ is even and
$\overline{\theta_I}=\overline{I}=1$ if $|I|$ is odd. Now we apply
Polya functor to obtain
\begin{prop}
The product rule in the algebra of symmetric odd functions
$(\bigwedge [\theta_1,\dots,\theta_m]^{\otimes n})_{S_n}$ is given
by
$$(\overline{\theta_{I_1}\dots \theta_{I_n}})(\overline{\theta_{J_1}\dots \theta_{J_n}})=\frac{1}{n!}\sum_{\sigma \in S_n}
\left( \sgn(I,J,\sigma) \prod_{k=1}^{n}
c(I_k,J_{\sigma^{-1}(k)})\right) \overline{\prod_{i=1}^{n}
\theta_{{I_k} \cup J_{\sigma^{-1}(k)}}},$$ where
$\sgn(I,J,\sigma)=(-1)^{e}$ and ${\des
e=\sum_{k>l}\overline{I_k}\overline{J_{\sigma^{-1}(l)}}+\sum_{\sigma
(k)>\sigma(l)} \overline{J_k J_l}}.$
\end{prop}
\begin{proof}
\begin{eqnarray*}
(\overline{\theta_{I_1}\dots
\theta_{I_n}})(\overline{\theta_{J_1}\dots
\theta_{J_n}})&=& \frac{1}{n!}\sum_{\sigma \in
S_n}\sgn(I,J,\sigma)(\theta_{I_1}\theta_{J_{\sigma^{-1}(1)}})\dots(\theta_{I_n}\theta_{J_{\sigma^{-1}(n)}})\\
\mbox{} &=& \frac{1}{n!}\sum_{\sigma \in
S_n}\sgn(I,J,\sigma)\left( \prod_{k=1}^{n}
c(I_k,J_{\sigma^{-1}(k)})\theta_{I_k\cup J_{\sigma^{-1}(k)}}\right) \\
\mbox{} &=& \frac{1}{n!}\sum_{\sigma \in S_n}\left(
\sgn(I,J,\sigma) \prod_{k=1}^{n} c(I_k,J_{\sigma^{-1}(k)})
\right)\overline{\prod_{k=1}^{n}\theta_{I_k\cup
J_{\sigma^{-1}(k)}}}
\end{eqnarray*}
\end{proof}

\subsection{Symmetric Boolean algebra}

Fix $n\in\mathbb{N}$ and let $P(n)$ be the free
$\mathbb{C}$-vector space generated by the subsets of $[n]$, i.e.,
$P(n)=\langle A: \ A\subset [n]\rangle$. Define a product $\cup$
on $P(n)$ by
$$\begin{array}{cccc}
  \cup: & P(n)\otimes P(n) & \longrightarrow & P(n) \\
   & A\otimes B & \longmapsto & A\cup B \\
\end{array}$$
$(P(n),\cup)$ is a Boolean algebra and $\dim(P(n))=2^{n}$. $S_n$
acts naturally on $[n]$ and thus on $P(n)$. We call the algebra
$(P(n),\cup)/S_n\cong (P(1)^{\otimes n})/S_n=\sym^{n}(P(1))$ the
{\em symmetric Boolean algebra}; it has dimension  $n+1$, a basis
being $\{[\overline{0}], [\overline{1}],\dots, [\overline{n}]\}$.
We define $P(n,k):=\{ A\subset [n]: |A|=k\}$. An application of
Polya functor yields the next

\begin{thm}

$[\overline{a}] [\overline{b}]= {\displaystyle \frac{1}{{n \choose
b}} \sum_{k=0}^{m} {a \choose b-k} {n-a \choose k}
[\overline{a+k}] }$, \ \ for all\ \
$[\overline{a}],[\overline{b}]\in (P(n),\cup)_{S_n}$, and\\
$m=\min(b,n-a)$.
\end{thm}
\begin{proof}
\begin{eqnarray*}
[\overline{a}] [ \overline{b}]&=& \frac{1}{n!}\sum_{\sigma \in
S_n} \overline{[a]\cup \sigma[b]} =\frac{1}{{n\choose
b}}\sum_{B\in P[n,b] } \overline{[a]\cup
B}\\
\mbox{} &\mbox{}& \mbox{}\\
\mbox{}&=& \frac{1}{{n\choose b}} \sum_{\begin{array}{c}
 {\scriptstyle B_0\subset P([n]-[a],k)} \\
 {\scriptstyle B_1\subset P([a],b-k)} \\
\end{array}}\overline{[a]\cup B_0}=\frac{1}{{n \choose b}} \sum_{k=0}^{m} {a \choose b-k}
{n-a \choose k} [\overline{a+k}]
\end{eqnarray*}
\end{proof}


\section{Quantum symmetric functions}

In this section we assume  the reader is familiar with the
notations from $\cite{Kon}$. Let us recall the notion of a formal
deformation
\begin{defi} Fix
a Poisson manifold $(M,\{ , \})$. A formal deformation
(deformation quantization) of the algebra of smooth functions on
$M$ is an associative star product\\ $\star: C^{\infty}(M)
[[\hbar]]\otimes_{\mathbb{R}[[\hbar]]} C^{\infty}(M)[[\hbar]]
\longrightarrow C^{\infty}(M)[[\hbar]]$ such that:
\begin{enumerate}
\item{ $f\star
g=\displaystyle{\sum_{n=0}^{\infty}B_n(f,g)\hbar^{n}}$, where
$B_n(-,-)$ are bi-differential operators .} \item{$f\star
g=fg+\frac{1}{2}\{f,g\}\hbar+O(\hbar^{2})$, where $O(\hbar^{2})$
are terms of order $\hbar^{2}$.}
\end{enumerate}
\end{defi}

In $\cite{Kon}$ a canonical $\star$-product has been constructed
for any Poisson manifold. For manifold $(\mathbb{R}^{m},\alpha)$
with Poisson bivector $\alpha$, the $\star$-product is given by
the formula
$$f\star g
=\sum_{n=0}^{\infty}\frac{\hbar^{n}}{n!}\sum_{\Gamma \in G_n}
\omega_\Gamma B_{\Gamma,\alpha}(f,g),$$ where $G_n$ is a
collection of {\em admissible graphs} each of which has $n$ edges,
and $\omega_{\Gamma}$ are some constants (independent of the
Poisson manifold). Given a finite group $K$ acting on
$(C^{\infty}(M),\star)$ by algebra automorphisms, we call the
algebra $(C^{\infty}(M)[[\hbar]] ,\star)_K\cong
(C^{\infty}(M)[[\hbar]],\star)^{K}$ the algebra of {\em quantum
$K$-symmetric functions} on $M$.

Next theorem shows how groups of automorphisms of
$(C^{\infty}(M)[[\hbar]],\star)$ arise in a natural way.

\begin{thm}
Assume we are given a Poisson structure $\{\mbox{-},\mbox{-}\}$ on
$\mathbb{R}^{m}$, and a group $K\subset S_m$ such that
$\{\mbox{-},\mbox{-}\}$ is $K$-equivariant. Then $K$ acts on
$(C^{\infty}(\mathbb{R}^{m})[[\hbar]],\star)$ by automorphisms.
\end{thm}
\begin{proof}
We assume that $\{f,g\}\circ \sigma=\{f\circ \sigma, g\circ
\sigma\}$, for all $f,g\in C^{\infty}(\mathbb{R}^{m})$, $\sigma\in
K$, or equivalently  $\alpha^{ij}(\sigma x)=\alpha^{\sigma(i)
\sigma(j)}$, where $\alpha^{ij}=\{x_i,x_j\}$ for all $i,j\in[m]$. Let us show that
$\sigma(f\star g)=(\sigma f)\star(\sigma g)$.

$$\sigma(f\star g)(x)=(f\star g)(\sigma^{-1} x)=\sum_{n=0}^{\infty}\frac{\hbar^{n}}{n!}
\sum_{\Gamma} w_{\Gamma} B_{\Gamma,\alpha}(f,g)(\sigma^{-1} x).$$
On the other hand
$$(\sigma f)\star(\sigma g)(x)=\sum_{n=0}^{\infty}\frac{\hbar^{n}}{n!} \sum_{\Gamma}
w_{\Gamma} B_{\Gamma,\alpha}(\sigma f,\sigma g)( x).$$ We need to
prove that
$$B_{\Gamma,\alpha}(f,g)(\sigma^{-1} x)=
B_{\Gamma,\alpha}(\sigma f,\sigma g)( x), \ \ \mbox{for all}\ \
\Gamma\in G_n.$$
Using Kontsevich's formula (see ${\cite{Kon}}$) we get
\begin{eqnarray*}
B_{\Gamma,\alpha}(f,g)(\sigma^{-1} x)&=&{\displaystyle
\sum_{I:E_{\Gamma} \longrightarrow [m]}
\left[\prod_{i=1}^{n}\left(\prod_{e\in E_{\Gamma}, e=(\ast,i)}
\partial _{I(e)}\right) \alpha^{I(e_i^{1})I(e_i^{2})}\right](\sigma^{-1}
x)\times}\\
\mbox{}& \mbox{}& \mbox{}\\
\mbox{}& \mbox{}&{\displaystyle \left( \left(\prod_{e\in E_{\Gamma}, e=(\ast,L)}
\partial _{I(e)}\right)   f\right)(\sigma^{-1} x)\times \left( \left(\prod_{e\in E_{\Gamma}, e=(\ast,R)}
\partial _{I(e)}\right)   g\right)(\sigma^{-1} x) }\\
\mbox{}& \mbox{}& \mbox{}\\
\mbox{}&=& {\displaystyle \sum_{I:E_{\Gamma} \longrightarrow [m]}
\left[\prod_{i=1}^{n}\left(\prod_{e\in E_{\Gamma}, e=(\ast,i)}
\partial _{\sigma(I(e))}\right) \alpha^{\sigma(I(e_i^{1}))\sigma(I(e_i^{2}))}\right](x)\times}\\
\mbox{}& \mbox{}& \mbox{}\\
\mbox{}& \mbox{}&{\displaystyle \left( \left(\prod_{e\in E_{\Gamma}, e=(\ast,L)}
\partial _{\sigma(I(e))}\right)   f\circ \sigma^{-1}\right)(x)\times \left( \left(\prod_{e\in E_{\Gamma}, e=(\ast,R)}
\partial _{\sigma(I(e))}\right)   g\circ \sigma^{-1} \right)(x) }\\
\mbox{}& \mbox{}& \mbox{}\\
\mbox{}&=& B_{\Gamma,\alpha}(\sigma f,\sigma g)(x).
 \end{eqnarray*}
\end{proof}

\begin{cor}
Under the conditions above, the product rule on
$(C^{\infty}((\mathbb{R}^{m})^{n})[[\hbar]],\star)_K$ is given by
\begin{equation}\label{ksp}
\overline{f}\star \overline{g}=\sum_{\sigma \in K}
\sum_{n=0}^{\infty} \frac{\hbar^{n}}{n!}\left( \sum_{\Gamma}
w_{\Gamma} \overline{B_{\Gamma, \alpha} (f,g\circ
\sigma^{-1})}\right)
\end{equation}
 for all $f,g\in
C^{\infty}(\mathbb{R}^{m})^{n}[[\hbar]]$.
\end{cor}

\begin{proof}Using Polya functor corollary \ref{prod2}, we have
$$ \overline{f}\star \overline{g}=\sum_{\sigma\in K} \overline{f\star \sigma g}=\sum_{\sigma \in K}
\sum_{n=0}^{\infty} \frac{\hbar^{n}}{n!}\left( \sum_{\Gamma}
w_{\Gamma} \overline{B_{\Gamma, \alpha} (f,g\circ
\sigma^{-1})}\right).$$

\end{proof}

\begin{defi}
Given a Poisson manifold $(\mathbb{R}^{m}, \{\mbox{ }, \mbox{
}\})$ and a subgroup $K\subset S_n$   the algebra of {\em quantum
symmetric functions} on $(\mathbb{R}^{m})^{n}$ is set to be
$(C^{\infty}(\mathbb{R}^{m})^{n}[[\hbar]],\star)_{K}\cong
(C^{\infty}(\mathbb{R}^{m})^{n}[[\hbar]],\star)^{K}$.
\end{defi}
Notice that if $(\mathbb{R}^{m},\alpha)$ is a Poisson manifold
then $(\mathbb{R}^{m})^{n}$ is a Poisson manifold in a natural
way. Moreover the Poisson structure on $(\mathbb{R}^{m})^{n}$ is
$S_n$-equivariant, and thus $K$-equivariant for all subgroup $K$
of $S_n$.
\subsection{Weyl algebra}
The Kontsevich $\star$-product given by formula ${\des f\star g
=\sum_{n=0}^{\infty}\frac{\hbar^{n}}{n!}\sum_{\Gamma \in G_n}
\omega_\Gamma B_{\Gamma,\alpha}(f,g)}$ is notoriously difficult to
compute. Nevertheless, there are two main examples, see \cite{Kon}
in which a fairly explicit knowledge of the start product is
available:
\begin{enumerate}
\item{If $\alpha$ is a constant non-degenerated Poisson bracket on $\mathbb{R}^{2n}$, then the quantum
algebra of polynomial functions on $\mathbb{R}^{2n}$, i.e.,
$(\mathbb{C}[x_1,\dots,x_{2n}][[\hbar]],\star)$ is isomorphic to $
W\otimes_{\mathbb{C}[[\hbar]]} \dots
\otimes_{\mathbb{C}[[\hbar]]} W$, where $W$ is the Weyl algebra, (see definition below). }
\item{If $\alpha$ is linear Poisson bracket in $\mathbb{R}^{n}$, then $(\mathbb{R}^{n},\alpha)$
is isomorphic as a Poisson manifold to $\mathfrak{g}^{\ast}$ for
some Lie algebra $\mathfrak{g}$. In this case the quantum algebra
of polynomial functions on $\mathfrak{g}^{\ast}$, i.e.,
$(\mathbb{C}[\mathfrak{g}^{\ast}][[\hbar]],\star)$ is isomorphic
to the universal enveloping algebra $U_h(\mathfrak{g})$ of
$\mathfrak{g}$.}

\end{enumerate}
Case $a$ will be considered in this section. Case $b$ for
$\mathfrak{g}=\mathfrak{sl}_2$ is considered in $\cite{DP1}$. The
case of a classical Lie algebra will be treated by our means
elsewhere. The algebra $W=\mathbb{C} \langle x,y \rangle[[\hbar]]/
\langle yx-xy-\hbar \rangle$ is called the {\em Weyl algebra}, it
is isomorphic to the canonical deformation quantization of
$(\mathbb{R}^{2},\dif x\wedge \dif y)$ if we consider only
polynomial functions on $\mathbb{R}^{2}$. This algebra admits a
natural representation as indicated in the
\begin{prop}\label{rweyl}
The map $\rho:W\longrightarrow
\en_{\mathbb{C}[[\hbar]]}(\mathbb{C}[x][[\hbar]])$ given by
$\rho(x)(f)=xf$ and $\rho(y)(f)=\hbar\frac{\partial f}{\partial
x}$, for any $f\in \mathbb{C}[x][[\hbar]]$  defines an irreducible
representation of the Weyl algebra.
\end{prop}
We order the letters of the Weyl algebra as follows: $x<y<\hbar$.
Assume we are given $A_i=(a_i,b_i)\in \mathbb{N}^{2}$, for
$i\in[n]$. Set $A=(A_1,\dots,A_n)\in (\mathbb{N}^{2})^{n}$,
$X^{A_{i}}=x^{a_i}y^{b_i}$ and let $|\mbox{
}|:\mathbb{N}^{n}\longrightarrow \mathbb{N}$ be the function such
that $|x|:=\sum_{i=1}^{n} x_i$, for all $x\in \mathbb{N}^{n}$.
Given  $x\in \mathbb{N}^{n}$ and $i\in\mathbb{N}$, we denote by
$x_{<i}$ the vector $(x_1,\dots,x_{i-1})\in \mathbb{N}^{i-1}$, by
$x_{\leq i}$ the vector $(x_1,\dots,x_{i})\in \mathbb{N}^{i}$ and
by $x_{>i}$ the vector $(x_{i+1},\dots,x_{n})\in
\mathbb{N}^{n-i}$. We write
$a\vdash n$ if $a\in\mathbb{N}^{k}$ for some $k$ and $|a|=n$.
Using this notation we have

\begin{defi}\label{ncweyl}
The {\em normal coordinates} $ N(A,k)$ of \ ${\des\prod_{i=1}^{n}
X^{A_i} \in W}$ are defined through the identity
\begin{equation}\label{nc}
\prod_{i=1}^{n} X^{A_i}=
 {\des \sum_{k=0}^{\min} N(A,k)x^{|a|-k}y^{|b|-k}\hbar^{k}}
\end{equation}
 for $0\leq k\leq \min=\min(|a|,
|b|)$. For $k>\min$, we set $N(A,k)$ equal to $0$.
\end{defi}
Recall that given finite sets $N$ and $M$ with $n$ and $m$
elements respectively, the number of one-to-one functions
$f:N\longrightarrow M$ is\ \ $m(m-1)\dots (m-n+1)=(m)_{n}$. The
number $(m)_{n}$ is called the $n$-th {\em falling factorial} of
$m$. The $n$-th {\em rising factorial} $m^{(n)}$  of $m$ is given
by $m^{(n)}=m(m+1)\dots (m+n-1)$. For any $a,b\in \mathbb{N}^{n}$,
we define ${a\choose b}:={a_1\choose b_1} {a_2\choose b_2}\dots
{a_{n}\choose b_{n}}$, and $a!=a_1!a_2!\dots a_n!$.

\begin{defi}
A k-pairing from set $E$ to set $F$ is an injective function from
a k-elements subset of $E$ to $F$. We denote by $P_{k}(E,F)$ the
set of k-pairings from $E$ to $F$.
\end{defi}

\begin{defi} Fix variables $t=(t_1,\dots, t_n)$ and $s=(s_1,\dots,
s_n)$. The {\em generating series} $N$ of the normal coordinates
in the Weyl algebra is given by
\begin{equation}\label{for}
N=\sum_{a,b,c} N(A,c) {\des\frac{s^{a}}{a!}\frac{t^{b}}{b!}
u^{c}}\in\mathbb{C}[[s,t,u]]
\end{equation}
where the sum runs over $a,b \in\mathbb{N}^{n}$, $c\in\mathbb{N}$
and $A=(A_1,\dots,A_n)\in(\mathbb{N}^{2})^{n}$.
\end{defi}

\begin{thm}\label{ecn}
Let $A,k$  be as in the Definition \ref{ncweyl}, the following
identity holds
\begin{enumerate}
\item{ $N(A,k)={\des\sum_{p\vdash k}{b\choose p}
\prod_{i=1}^{n-1}\left(|a_{>i}|-|p_{>i}| \right)_{p_i}}$, where
$p\in \mathbb{N}^{n-1}$.}
\item{ Let
$E_1,\dots,E_n,F_1,\dots,F_n$ be disjoint sets such that
$\sharp(E_i)=a_i$, $\sharp(F_i)=b_i$, for $i\in[n]$, set $
E=\cup_{i=1}^{n} E_i$, and $F=\cup_{i=1}^{n} F_i$, then $ N(A,k)=
\sharp(\{ p \in P_{k}(E,F) \ | \ \mbox{if} \ \ (a,p(a)) \in {E_i
\times F_j} \ \mbox{then} \ i>j \}).$}
\item{$N = \exp\left(
{\des\sum_{i>j}ut_is_j+\sum_i t_i+\sum_j s_j }\right).$}
\end{enumerate}
\end{thm}
\begin{proof} Using induction one show that the following identity
hold in the Weyl algebra
\begin{equation}\label{bas}
y^{b}x^{a}=\sum_{k=0}^{\min}{b\choose k}(a)_k x^{a-k}
y^{b-k}\hbar^{k},
\end{equation} where $\min=\min(a,b)$ . Several applications of
identity $(\ref{bas})$ imply $a.$ Notice that for given sets $E,F$
such that $\sharp(E)=a$ and $\sharp(F)=b$, ${b\choose k}(a)_k$ is
equal to $\sharp(\{p\in P_k(E,F)\})$, showing part $b$ for $n=2$.
The general formula follows from induction. This prove $b.$ It
follows from standard combinatorial facts (see $\cite{RS}$,
$\cite{KW}$ for more details)  that
$$\exp\left( {\des\sum_{i>j}ut_is_j+\sum_i t_i+\sum_j s_j
}\right)=\sum_{a,b,c} c_{a,b,c}\frac{s^{a}}{a!}\frac{t^{b}}{b!}
\frac{u^{c}}{c!}$$ where $c_{a,b,c}=\sharp(\{(p,\sigma): p\in
P_c(A,B) \ \mbox{and}\ \sigma:[c]\rightarrow p, \ \mbox{a
bijection}\ \})$ which is equivalent to formula $(\ref{for})$.
\end{proof}

\begin{figure}[ht]
\begin{center}
\includegraphics[width=3.5in]{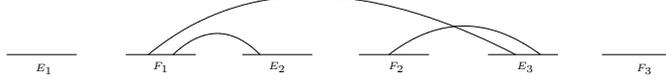}
\caption{Combinatorial interpretation of N. \label{fig:weyl}}
\end{center}
\end{figure}

Figure $\ref{fig:weyl}$ illustrates the combinatorial
interpretation of the normal coordinate $N(A,k)$ of an element
\\ ${\des\prod_{i=1}^{3}X^{A_i}\in W} $, it shows are of the
possible pairing contributing to $N(A,k)$.

\begin{cor} For any given $(t,a,b)\in \mathbb{N}\times\mathbb{N}^{n}\times \mathbb{N}^{n}$, the
following identity holds $$\prod_{i=1}^{n} \left(t+|a_{>i}|-
|b_{>i}|\right)_{b_i}=\sum_{p\vdash k}{b\choose p}
\prod_{i=1}^{n-1}\left(|a_{>i}|-|p_{>i}| \right)_{p_i}t_{|b|
-k}$$
\end{cor}
\begin{proof} Consider the identity $(\ref{nc})$ in the representation of
the Weyl algebra defined in Proposition $\ref{rweyl}$. Apply both
sides of the identity $(\ref{nc})$ to $x^{t}$ for $t\in
\mathbb{N}$ and use Theorem $\ref{ecn}$.
\end{proof}

\begin{thm}
Let $A,c$  be as in the Definition \ref{ncweyl}. The following
identity holds
$$ N(A,c)={\des\sum_{\sum
c_{ij}=c} \prod_{j=2}^{k} { a_j \choose \sum_i c_{ij}}{ \sum_i
c_{ij} \choose c_{i j}\dots  c_{(j-1) j}} \prod_{i=1}^{k-1} { b_i
\choose \sum_j c_{ij}}{ \sum_j c_{ij} \choose c_{i (i+1)}\dots
c_{i k}}\prod_{i>j} c_{ij}! } $$ where $c_{ij}$ are integers with
$1\leq i\leq k-1$, $2\leq j\leq k$ and $i>j$.
\end{thm}
\begin{proof}
Consider maps $F_{ij}:\{ p \in P_{c}(E,F) \ | \ \mbox{if} \ \
(a,p(a))
\in {E_i
\times F_j} \ \mbox{then} \ i>j \}\longrightarrow \mathbb{N}$
given by $F_{ij}(p)=\sharp (\{(a,p(a))\in {E_i\times F_j}\})$.
Notice that $N(A,c)={\des \sum_{\sum{c_{ij}=c}} \sharp(\{ p\in
P_c(E,F): F_{ij}(p)=c_{ij}, i>j\})}$. Moreover
$$\sharp(\{ p\in
P_c(E,F): F_{ij}(p)=c_{ij}, i>j\})=\!\!\prod_{j=2}^{k} { a_j
\choose
\sum_i c_{ij}}{ \sum_i c_{ij} \choose c_{i j}\dots  c_{(j-1) j}}
\prod_{i=1}^{k-1} { b_i
\choose \sum_j c_{ij}}{ \sum_j c_{ij} \choose c_{i (i+1)}\dots
c_{i k}}\prod_{i>j} c_{ij}!. $$
\end{proof}

\subsection{Quantum symmetric functions of Weyl type}\label{fstw}

The following theorem provides explicit formula for the product of
$m$ elements of the $n$-th symmetric power of the Weyl algebra.
Let us explain our notation: fix a matrix $A:[m]\!\times\!
[n]\longrightarrow \mathbb{N}^{2}$,
$(A_{ij})=((a_{ij}),(b_{ij}))$. Given $\sigma \in (S_n)^{m}$ and
$j\in[n]$, $A_j^{\sigma}$ denotes the vector $(A_{1
\sigma_{1}^{-1}(j)},\dots,A_{m\sigma_{m}^{-1}(j)})\in
(\mathbb{N}^{2})^{n}$ and set $X_j^{A_{ij}}=x_j^{a_{ij}}
y_j^{b_{ij}}$ for $j\in[n]$. Set $|A_j^{\sigma}|=(|a_j^{\sigma}|,
|b_j^{\sigma}|)$ where ${\des |a_j^{\sigma}|=\sum_{i=1}^{m}
a_{i\sigma^{-1}_i(j)}}$ and ${\des |b_j^{\sigma}|=\sum_{i=1}^{m}
b_{i\sigma^{-1}_i(j)}}$. We have the following

\begin{thm}\label{sp}For any  $A:[m]\!\times\!
[n]\longrightarrow \mathbb{N}^{2}$, the following identity
\begin{equation}
{\des (n!)^{m-1}\prod_{i=1}^{m}\left(\overline{\prod_{j=1}^{n}
X_j^{A_{ij}}}\right)= \sum_{\sigma,k,p}\left(\prod_{i,j} {
b_j^{\sigma} \choose p^{j} }
({|(a_{j}^{\sigma})}_{>i}|-|p^{j}_{>i}|)_{p_i^{j}}\right)
\overline{\prod_{j=1}^{n} X_j^{|A_j^{\sigma}|-(k_j,k_j)}}
\hbar^{|k|}}
\end{equation} where
$\sigma\in {\{\id\} \times S_n^{m-1}}$, $k\in\mathbb{N}^{n}$,
$(i,j)\in [m-1]\times [n]$ and $p=p_i^{j}\in
(\mathbb{N}^{m-1})^{n}$, holds in $\sym^{n}(W)$.
\end{thm}
\begin{proof} We use Theorem \ref{pfu} and Theorem \ref{ecn}
\begin{eqnarray*}
{\des (n!)^{m-1}\overline{\prod_{i=1}^{m}\left(\prod_{j=1}^{n}
X_j^{A_{ij}}\right)}}&=&\sum_{\sigma \in \{\id\} \times
S_n^{m-1}}\overline{\prod_{j=1}^{n}\left( \prod_{i=1}^{m}
X_j^{A_{i
\sigma^{-1}_i(j)}}\right)}\\
\mbox{}&=&\sum_{\sigma \in \{\id\} \times
S_n^{m-1}}\overline{\prod_{j=1}^{n}\left(\sum_{k=0}^{\min_j}
N(A_j^{\sigma},k)X_j^{|A_j^{\sigma}|-(k,k)} \hbar^{k}  \right)}\\
\mbox{}&=&\sum_{\sigma,k}\left(\prod_{j=1}^{n}N(A_j^{\sigma}
,k_j)  \right)  \overline{\prod_{j=1}^{n}
X_j^{|A_j^{\sigma}|-(k_j,k_j)}}
\hbar^{|k|}\\
\mbox{}&=& \sum_{\sigma,k,p}\left(\prod_{i,j} { b_j^{\sigma}
\choose p^{j} }
({|(a_{j}^{\sigma})}_{>i}|-|p^{j}_{>i}|)_{p_i^{j}}\right)
\overline{\prod_{j=1}^{n} X_j^{|A_j^{\sigma}|-(k_j,k_j)}}
\hbar^{|k|},
\end{eqnarray*}
where $\min_j=\min(|a_j^{\sigma}|,|b_j^{\sigma}|)$.
\end{proof}

Now we proceed to state and prove the quantum analogue of Theorem
\ref{csf}. We begin by introducing a $\star$-product on
$(\mathbb{C}[\mathbb{N}^{2m}]^{\otimes n}[[\hbar]])_{S_n}$, which
is motivated by the proof of
 Theorem \ref{tt} below.

\begin{defi}
The $\star$-product on $(\mathbb{C}[\mathbb{N}^{2m}]^{\otimes
n}[[\hbar]])_{S_n}$ for $A\in (\mathbb{N}^{2m})^{n}$ and $C\in
(\mathbb{N}^{2m})^{n}$ is given by the formula

\begin{equation}\label{pqs}
\overline{A }\star\overline{C}={\displaystyle
\frac{1}{n!}\sum_{I,\sigma} {b \choose I} (\sigma(c))_{I}
\overline{{A+\sigma(C)-(I,I)}} \hbar^{|I|}}\end{equation} where
$I:[n]\times [m]\rightarrow \mathbb{N}$ and $\sigma\in S_n$.
\end{defi}

\begin{defi} Fix $m\in\mathbb{N}$. The algebra of quantum symmetric functions
of type $A_n$, $B_n$ and $D_n$ are given by
\begin{itemize}
\item{$\qsym_{A_n}(m)=(\mathbb{C}[x_1,y_1,\dots,x_n,y_n][[\hbar]],\star)_{S_n}\cong
(\mathbb{C}[x_1,y_1,\dots,x_n,y_n][[\hbar]],\star)^{{S_n}}$,}
\item{$\qsym_{B_n}(m)=
(\mathbb{C}[x_1,y_1,\dots,x_n,y_n][[\hbar]],\star)_{{\mathbb{Z}_2^{n}\rtimes
S_n}}\cong
(\mathbb{C}[x_1,y_1,\dots,x_n,y_n][[\hbar]],\star)^{{{\mathbb{Z}_2^{n}\rtimes
S_n}}}$,}
\item{$\qsym_{D_n}(m)=
(\mathbb{C}[x_1,y_1,\dots,x_n,y_n][[\hbar]],\star)_{{\mathbb{Z}_2^{n-1}\rtimes
S_n}}\cong
(\mathbb{C}[x_1,y_1,\dots,x_n,y_n][[\hbar]],\star)^{{{\mathbb{Z}_2^{n-1}\rtimes
S_n}}}$,}
\end{itemize}
where $x_i=(x_{i1},\dots,x_{im})$ and $y_i=(y_{i1},\dots,y_{im})$.
\end{defi}

\vspace{0.3cm}

\begin{thm}[Quantum symmetric functions]\label{tt}\quad
\begin{enumerate}
\item{The following is a commutative diagram
$$\begin{array}{ccccc}
  (\mathbb{C}[\mathbb{N}^{2m}]^{\otimes n}[[\hbar]],\star)_{S_n} & \supset & (\mathbb{C}[\mathbb{N}^{2m}_e +
\mathbb{N}^{2m}_o]^{\otimes n}[[\hbar]],\star)_{S_n} & \supset & (\mathbb{C}[\mathbb{N}^{2m}_e]^{\otimes n}[[\hbar]],\star)_{S_n} \\
  \downarrow &  & \downarrow &  & \downarrow \\
  \qsym_{A_n} (m) &\twoheadrightarrow &  \qsym_{D_n}(m) & \twoheadrightarrow & \qsym_{B_n}(m) \\
\end{array}$$ where the vertical arrows are isomorphisms.}

\item{For $A_i=(a_i,b_i)\in (\mathbb{N}^{2})^{m}$
and $C_i=(c_i,d_i)\in (\mathbb{N}^{2})^{m}$ $i\in[n]$, we set
$X_i^{A_i}=x_i^{a_i}y_i^{b_i}$ and $X_i^{C_i}=x_i^{c_i}y_i^{d_i}$
where $x_i=(x_{i1},\dots,x_{im})$ and $y=(y_{i1},\dots,y_{im})$,
the product rule in $\qsym_{A_n}(m)$ is given by
 $$\overline{X_1^{A_1} X_2^{A_2}\dots X_n^{A_n}}\star  \overline{X_1^{C_1}
X_2^{C_2}\dots X_n^{C_n}}={\displaystyle
\frac{1}{n!}\sum_{I, \sigma } {b \choose I} (\sigma(c))_{I}  \overline{X^{A+\sigma(C)-(I,I)}} \hbar^{|I|}},$$\\
where $I:[n]\times [m]\rightarrow \mathbb{N}$, $\sigma\in S_n$.
The product on $\qsym_{B_n}(m)$ and \ $\qsym_{D_n}(m)$ is given
by restriction.}
\end{enumerate}
\end{thm}
\begin{proof}

We prove $b$ which implies $a$. Given $A_i,C_i\in
(\mathbb{N}^{2})^{m}$, using $(\ref{pqs})$, we obtain

\begin{eqnarray*}
\overline{X_1^{A_1} X_2^{A_2}\dots X_n^{A_n}} \star
\overline{X_1^{C_1} X_2^{C_2}\dots X_n^{C_n}}&=&
\frac{1}{n!}\sum_{\sigma\in S_n} \overline{X_1^{A_1} X_1^{C_{ \sigma^{-1}
(1)}} X_2^{A_2} X_2^{C_{ \sigma^{-1} (2)}}\dots X_n^{A_n} X_n^{C_{
\sigma^{-1}
(n)}}}\\
\mbox{}&\mbox{}&\mbox {}\\
\mbox{} &=& \frac{1}{n!}\sum_{\sigma \in S_n} \overline{x_1^{a_1} y_1^{b_1}  x_1^{c_{\sigma^{-1}(1)}}
 y_1^{d_{\sigma ^{-1}(1)}}\dots  x_n^{a_n} y_n^{b_n}
x_n^{c_{\sigma^{-1}(n)}} y_n^{d_{\sigma^{-1} (n)}}}\\
\mbox{}&\mbox{}&\mbox {}\\
\mbox{} &=&\frac{1}{n!}\sum_{\sigma \in S_n} \left( {\des
\prod_{j=1}^{n} \sum_{i_j=0}^{\min_j} {b_j \choose i_j}
({c_{\sigma^{-1}(j)}})_{i_j}
 \overline{X_j^{A_j+C_{\sigma^{-1}(j)}-(i_j,i_j)}}\hbar^{i_j}} \right)\\
\mbox{}&\mbox{}&\mbox {}\\
\mbox{} &=&\frac{1}{n!}\sum_{I, \sigma} \left( {b \choose
I}(\sigma(c))_I \right) \overline{X^{A+\sigma(C)-(I,I)}}
\hbar^{|I|}
\end{eqnarray*}
where $\min_j=\min(b_j,c_{\sigma^{-1}(j)})$. Now, consider
$A_i,C_i \in \mathbb{N}^{2m}_e$

\begin{eqnarray*}
\overline{X_1^{A_1}\dots X_n^{A_n}}\star \overline{X_1^{C_1} \dots
X_n^{C_n}}&=&\frac{1}{2^{n} n!}\sum_{(t,\sigma) \in
\mathbb{Z}_2^{n }\rtimes S_n} \overline{(x_1^{a_1}y_1^{b_1} \dots
x_n^{a_n}y_n^{b_n}) ((t,\sigma)x_1^{c_1}y_1^{d_1} \dots
x_n^{c_n}y_n^{d_n})}\\
\mbox{}&\mbox{}&\mbox {}\\
\mbox{}&=&{\displaystyle \frac{1}{2^{n}n!} \sum_{I, \sigma }\left(
\sum_{t\in \mathbb{Z}_2^{n}} (-1)^{{\displaystyle\sum_{t_i=-1}
c_{i}+d_{i}}}{b \choose I}(\sigma(c))_{I} \right)
\overline{X^{A+\sigma(C)-(I,I)}}
\hbar^{|I|} }\\
\mbox{}&\mbox{}&\mbox {}\\
\mbox{}&=&{\displaystyle \frac{1}{n!}\sum_{I, \sigma }\left( {b
\choose I}(\sigma(c))_{I}\right)  \overline{X^{A+\sigma(C)-(I,I)}}
\hbar^{|I|}}
\end{eqnarray*}
since \  $\sum_{t\in \mathbb{Z}_2^{n}}
(-1)^{{\displaystyle\sum_{t_i=-1} c_{i}+d_{i}}}=2^{n}$. Finally,
consider $A_i,C_i \in \mathbb{N}^{2m}_e +
\mathbb{N}^{2m}_o$

\begin{eqnarray*}
\overline{X_1^{A_1}\dots X_n^{A_n}}\star  \overline{X_1^{C_1} \dots
X_n^{C_n}}&=&\frac{1}{2^{n-1} n!}\sum_{(t,\sigma) \in
\mathbb{Z}_2^{n }\rtimes S_n} \overline{(x_1^{a_1}y_1^{b_1} \dots
x_n^{a_n}y_n^{b_n}) (t,\sigma)x_1^{c_1}y_1^{d_1} \dots
x_n^{c_n}y_n^{d_n}}\\
\mbox{}&\mbox{}&\mbox {}\\
\mbox{}&=& \frac{1}{2^{n-1}n!} \sum_{I, \sigma }\left( k(t,c,d){b
\choose I}(\sigma(c))_{I} \right) \overline{X^{A+\sigma(C)-(I,I)}}
\hbar^{|I|}\\
\mbox{}&\mbox{}&\mbox {}\\
\mbox{}&=&{\displaystyle \frac{1}{n!}\sum_{I, \sigma }\left( {b
\choose I}(\sigma(c))_{I}\right)  \overline{X^{A+\sigma(C)-(I,I)}}
\hbar^{|I|}}
\end{eqnarray*}
since \  $k(t,c,d)=\des\sum_{\begin{array}{c}
  {\scriptstyle t\in \mathbb{Z}_2^{n}} \\
  {\scriptstyle \prod t_i=1} \\
\end{array}} (-1)^{\sum_{t_i=-1}
c_{i}+d_{i}+c_n+d_n}=2^{n-1}$ \end{proof}

\subsection{Quantum symmetric functions on  $\mathbb{C}^{n}/\mathbb{Z}_m^{n}\rtimes S_n$}

In this section we shift to complex analytic notation to study the
Poisson orbifolds $\mathbb{C}^{n}/\mathbb{Z}_m^{n}\rtimes S_n$ and
$\mathbb{C}^{n}/ \mathcal{D}_m^{n}\rtimes S_n$. The deformation
quantization of $\mathbb{C}^{n}$ provided with the canonical
symplectic structure is isomorphic to the Weyl  algebra
$W=\mathbb{C}\langle z_1,
\overline{z_1},\dots ,z_n, \overline{z_n} \rangle[[\hbar]]/\langle z
\overline{z}-\overline{z} z-2i \hbar \rangle$.
The group $\mathbb{Z}_m^{n}\rtimes S_n$ acts on $\mathbb{C}^{n}$
as follows
$$\begin{array}{cccc}
 ( \mathbb{Z}_m^{n}\rtimes S_n)\times  \mathbb{C}^{n}& \longrightarrow & \mathbb{C}^{n}   \\
  ((w_1, \dots, w_n), \sigma)(z_1,\dots, z_n) & \longmapsto & (w_1z_{\sigma^{-1}(1)},\dots, w_n z_{\sigma^{-1} (n)}) \\
\end{array}$$
where $w_j=e^{\frac{2 \pi i k_j}{m}}$, $j=1,\dots,n$. Thus
$\mathbb{Z}_m^{n}\rtimes S_n$ acts on $W=\mathbb{C}\langle z_1,
\overline{z_1},\dots ,z_n, \overline{z_n} \rangle[[\hbar]]/\langle z
\overline{z}-\overline{z} z-2i \hbar \rangle$.\\ We denote
$Z_i^{A_i}= z_i^{a_i} {\overline{z}}_i^{b_i}$ and
$\mathbb{N}_m^{2}=\{(a,b): \ \mbox{there is}\ k\in\mathbb{Z} \
\mbox{such that}\ b-a=km \}.$

\begin{defi} The $\star$-product on $(\mathbb{C}[
\mathbb{N}_m^{2}]^{\otimes n}[[\hbar]])_{\mathbb{Z}_m^{n}\rtimes
S_n}$ for $A\in (\mathbb{N}^{2})^{n}$ and $C\in
(\mathbb{N}^{2})^{n}$ is given by the formula
$$\overline{A }\star\overline{C}={\displaystyle
\frac{1}{n!}\sum_{I,\sigma}(-2i)^{|I|} {b \choose I}
(\sigma(c))_{I} \overline{{A+\sigma(C)-(I,I)}} \hbar^{|I|}}$$
where  $I:[n] \rightarrow \mathbb{N}$ and $\sigma\in S_n$.
\end{defi}

\begin{thm} The map
$$\begin{array}{ccc}
 (\mathbb{C}[ \mathbb{N}_m^{2}]^{\otimes
n}[[\hbar]],\star)_{\mathbb{Z}_m^{n}\rtimes S_n}  &
\longrightarrow & ({\mathbb{C}} [z,\overline{z}]^{\otimes
n}[[\hbar]],\star)_{\mathbb{Z}_m^{n}\rtimes S_n}  \\
  \mbox{} & \mbox{} & \mbox{}  \\
  (A_1,\dots,A_n) & \longmapsto & Z_1^{A_1}\dots Z_n^{A_n}  \\
\end{array}$$
is an algebra isomorphism.
\end{thm}
\begin{proof} Let $A_i=(a_i,b_i)\in \mathbb{N}^{2}$, and
$C_{i}=(c_i,d_i)\in \mathbb{N}^{2}$ $i\in[n]$, we have
\begin{eqnarray*}
\overline{Z_1^{A_1}\dots Z_n^{A_n}} \star  \overline{Z_1^{C_1}\dots
Z_n^{C_n}}&=& \overline{z_1^{a_1} \overline{z}_1^{b_1} \dots
z_n^{a_n} \overline{z}_n^{b_n}} \  \overline{ z_1^{c_1}
\overline{z}_1^{d_1} \dots z_n^{c_n}
\overline{z}_n^{d_n}}\\
\mbox{} &\mbox{}& \mbox{}\\
\mbox{}&=&\frac{1}{m^{n} n!}\sum_{(w,\sigma)\in \mathbb{Z}_m^{n}\rtimes S_n}
 \overline{z_1^{a_1} \overline{z}_1^{b_1} \dots z_n^{a_n} \overline{z}_n^{b_n}
((w,\sigma)( z_1^{c_1} \overline{z}_1^{d_1} \dots z_n^{c_n}
\overline{z}_n^{d_n}))}\\
\mbox{} &\mbox{}&\mbox{}\\
\mbox{} &=& \frac{1}{m^{n} n!}\sum_{\sigma\in S_n}  \left( \prod_{j=1}^{n}\sum_{k_j=0}^{m-1}
(e^{\frac{2 \pi i }{m}(d_j-c_j)})^{k_j} \right)
 \overline{\prod_{j=1}^{n} z_j^{a_j} \overline{z}_j^{b_j} z_j^{c_{\sigma^{-1}(j)}}
 \overline{z}_j^{d_{\sigma^{-1} (j)}}}\\
 \mbox{} &\mbox{}&\mbox{}\\
 \mbox{}&=&\frac{1}{ n!}\sum_{\begin{array}{c}
  {\scriptstyle I:[n]\rightarrow
\mathbb{N}} \\
{\scriptstyle \sigma \in S_n} \\
\end{array}}\left(
(-2i)^{|I|}{b \choose I} {(\sigma(c))_{I}} \right)
\overline{Z^{A+\sigma(C)-(I,I)}} \hbar^{|I|}
\end{eqnarray*}
\end{proof}
Let $\mathcal{D}_m$ be the dihedral group of order $2n$, where
$\mathcal{D}_m=\{R_0,R_1,\dots,R_{m-1},S_0,S_1,\dots,S_{m-1}\}$,
$R_k(z)=e^{\frac{2\pi i k}{n}}z$, and $S_k(z)=e^{\frac{2\pi i
k}{n}}\overline{z},$ for $k=0,\dots,m-1$.
$\mathcal{D}_m^{n}\rtimes S_n$ acts on $\mathbb{C}^{n}$ as follows
$$\begin{array}{ccc}
  (\mathcal{D}_m^{n}\rtimes S_n)\times\mathbb{C}^{n}  & \longrightarrow & \mathbb{C}^{n} \\
  ((d_1,d_2,\dots,d_n),\sigma)(z_1,z_2,\dots,z_n) & \longmapsto &
  (d_1 z_{\sigma^{-1}(1)},\dots, d_n z_{\sigma^{-1}(n)})
\end{array}  $$
We will apply Polya functor to obtain a product rule on
$(C^{\infty}(\mathbb{C}^{n})[[\hbar]],\star)_{\mathcal{D}^{n}\rtimes
S_n}$

\begin{defi} The $\star$-product on $(\mathbb{C}[
\mathbb{N}_m^{2}]^{\otimes
n}[[\hbar]],\star)_{\mathcal{D}^{n}\rtimes S_n}$ for $A\in
(\mathbb{N}^{2})^{n}$ and $C\in (\mathbb{N}^{2})^{n}$ is given by
the formula
$$\overline{A }\star \overline{C}={\des \frac{1}{2 n!}\sum_{I,\sigma}{b\choose I}{(\sigma(c))_{I}} \overline{{A+(\sigma(c),
\sigma(d))-(I,I)}} \hbar^{|I|}} + {\des \frac{1}{2
n!}\sum_{I,\sigma}{b+c\choose I}{(\sigma(d))_{I}}
\overline{{A+(\sigma(d),\sigma(c))-(I,I)}}} \hbar^{|I|}$$ where
$I:[n] \rightarrow \mathbb{N}$ and $\sigma\in S_n$.
\end{defi}

\begin{thm} The map $(\mathbb{C}[
\mathbb{N}_m^{2}]^{\otimes
n}[[\hbar]],\star)_{\mathcal{D}^{n}\rtimes S_n}\longrightarrow
({\mathbb{C}} [z,\overline{z}]^{\otimes
n}[[\hbar]],\star)_{\mathcal{D}^{n}\rtimes S_n}$ given by\\ $
(A_1,\dots,A_n)  \longmapsto  Z_1^{A_1}\dots Z_n^{A_n}$ is an
algebra isomorphism.
\end{thm}

\begin{proof}  For any $A_i=(a_i,b_i)\in\mathbb{N}^{2}$,
$C_i=(c_i,d_i)\in\mathbb{N}^{2}$, we have\\

 $\overline{Z_1^{A_1}\dots
Z_n^{A_n}} \ \overline{Z_1^{C_1}\dots Z_n^{C_n}}$
\begin{eqnarray*}
\mbox{}&=&
\frac{1}{(2n) n!}\sum_{(d,\sigma)\in \mathcal{D}^{n}\rtimes S_n
}\overline{z_1^{a_1} \overline{z}_1^{b_1}
 \dots z_n^{a_n} \overline{z}_n^{b_n}((d,\sigma)  z_1^{c_1} \overline{z}_1^{d_1} \dots z_n^{c_n}
\overline{z}_n^{d_n})}\\
\mbox{}&\mbox{}&\mbox{}\\
\mbox{}&=&\frac{1}{(2n) n!}\sum_{\sigma\in S_n}\left( \prod_{j=1}^{n}\sum_{k_j=0}^{n-1}
e^{\frac{2 \pi i k_j}{n}(d_j-c_j)}\right) \left(
\prod_{j=1}^{n} \overline{z_j^{a_j} \overline{z_j}^{b_j}z_j^{c_{\sigma^{-1}(j)}}\overline{z_j}^{d_{\sigma^{-1}(j)}}}
+\overline{z_j^{a_j}
\overline{z_j}^{b_j+c_{\sigma^{-1}(j)}}z_j^{d_{\sigma^{-1}(j)}}}\right)\\
\mbox{}&\mbox{}&\mbox{}\\
\mbox{}&=&{\des \frac{1}{2 n!}\sum_{\sigma, I}{b\choose I}{(\sigma(c))_{I}} \overline{Z^{A+(\sigma(c),\sigma(d))-(I,I)}} \hbar^{|I|}} +
{\des \frac{1}{2 n!}\sum_{\sigma, I}{b+c\choose
I}{(\sigma(d))_{I}}
\overline{Z^{A+(\sigma(d),\sigma(c))-(I,I)}}
\hbar^{|I|}}.
\end{eqnarray*}
where $\sigma \in S_n$ and $I:[n]\rightarrow
\mathbb{N}$.
\end{proof}

\subsection{Quantum super-functions}\label{qspf}

We denote by $\mat(n)$ the algebra of $\mathbb{C}$-matrices of
order $n\times n$.  It is well-known see \cite{Del} that the
canonical quantization of
$\bigwedge_{\mathbb{C}}[\theta_1,\dots,\theta_m]$ is isomorphic to
the Clifford algebra $C(m)$, i.e, the complex free algebra on $m$
generators $\theta_1,\dots,\theta_m$ subject to the relations
$\theta_i \theta_j+ \theta_j\theta_i=2\delta_{ij}$, for
$i,j\in[m]$. It is also known that
$$\begin{array}{ll}
C(m)\cong&\mat(2^{\frac{m}{2}}) \ \mbox{if}\ m \ \mbox{is
even, and}\\
\mbox{}& \mbox{}\\
C(m) \cong&\mat(2^{\frac{m-1}{2}})\oplus \mat(2^{\frac{m-1}{2}}) \
\mbox{if}\ m \ \mbox{is odd.}
\end{array}$$
Thus the algebra of quantum symmetric functions on the super-space
$\mathbb{R}^{0|m}$ is isomorphic to
\begin{eqnarray*}
\sym^{n}(\mat(2^{m}))&\cong&\schur (n,2^{m})\ \mbox{for}\ m \ \mbox{even,}\\
{\des\sym^{n}(\mat(2^{\frac{m-1}{2}})\oplus
\mat(2^{\frac{m-1}{2}}))}&\cong &{\des\bigoplus_{a+b=n}
\sym^{a}(\mat(2^{\frac{m-1}{2}}))\otimes
\sym^{b}(\mat (2^{\frac{m-1}{2}}))}\\
\mbox{}&\cong&{\des\bigoplus_{a+b=n}\schur(a,2^{\frac{m-1}{2}})\otimes
\schur(b,2^{\frac{m-1}{2}})}\ \mbox{for}\ m \ \mbox{odd}.
\end{eqnarray*}
Definition \ref{glin} and Proposition \ref{glin1} below are taken
from $\cite{Mer}$.
\begin{defi}\label{glin}
We denote by $\mathfrak{gl}(\infty)$ the algebra of all matrices
$(a_{ij})$ such that $a_{ij}\in \mathbb{C}[[\hbar]]$ if $i\geq j$,
and $a_{ij}\in (\hbar \mathbb{C}[[\hbar]])^{j-i}$, if $i<j$. We
set $\schur(\infty,n):=\sym^{n}(\mathfrak{gl}(\infty))$.
\end{defi}
For any $a,b\in \mathbb{Z}^{\geq 0}$, we define the matrices
$$E_{a,b}(i,j)=\left \{ \begin{array}{cc}
  \frac{(b+k)!}{k!} \hbar^{b}, & \mbox{if}\ (i,j)=(a,b)+(k,k),\ k=0,1,2,\dots \\
  \mbox{} & \mbox{}\\
  0, & \mbox{otherwise}. \\
\end{array} \right.$$
\begin{prop}\label{glin1}
There is a canonical isomorphism
  $\rho:  W
  \longrightarrow  \mathfrak{gl}(\infty) $
defined on generators by $\rho(x)=E_{1,0}$, $\rho(y)=E_{0,1}$ and
$\rho(\hbar)=\hbar I$.

\end{prop}
\begin{proof}
The linear map $\rho$  is a well defined algebra homomorphism
since the Weyl algebra is the quotient of a (formal) free algebra
by the ideal generated by the relation $yx=xy+\hbar$, and the
following identity holds $E_{0,1} E_{1,0}=E_{1,0}E_{0,1}+\hbar I$
in $\mathfrak{gl}(\infty)$. Similarly $\rho$ is a bijection since
$E_{a,b}$ is a basis for $\mathfrak{gl}(\infty)$ and
$\rho(x^{a}y^{b})=E_{a,b}$ as consequence of the fact that in
$\mathfrak{gl}(\infty)$ the following identities are satisfied:
$E_{0,a}E_{0,1}=E_{0,a+1}$, $E_{a,0}E_{1,0}=E_{a+1,0}$ and
$E_{a,0}E_{0,b}=E_{a,b}$, for all $a,b\in\mathbb{Z}^{\geq 0}$.

\end{proof}
We denote by $\widehat{W}$ the formal Weyl algebra, i.e.
$\widehat{W}=\mathbb{C}\langle \langle x,y \rangle
\rangle[[\hbar]]/\langle yx-xy-\hbar\rangle$ and
$\hat{\mbox{\rm{Q}}}\sym_{A_n}(1)=\sym^{n}(\widehat{W})$.

\begin{thm}
\begin{enumerate}
\item{The algebra  $ \hat{\mbox{\rm{Q}}}\sym_{A_n}(1)$  of formal
quantum symmetric functions on \\ $((\mathbb{R}^{2})^{n},\sum \dif
x_i\wedge \dif y_i)$ is isomorphic to $\schur(\infty,n)$.}
\item{The algebra of formal quantum symmetric functions on the
superspaces $\mathbb{R}^{2|n}$ is isomorphic to
$\sym^{n}(\mathfrak{gl}(\infty)\otimes C(m))$.}
\end{enumerate}
\end{thm}
\begin{proof}

\begin{enumerate}
\item{$\hat{\mbox{\rm{Q}}}\sym_{A_n}(1)\cong \sym^{n}(\widehat{W})\cong \schur (\infty,n)$.}
\item{ The algebra of formal quantum functions on the superspace
$\mathbb{R}^{2|n}$ is isomorphic to
$$\widehat{W}\otimes C(m) \cong \mathfrak{gl}(\infty)\otimes C(m)
.$$ Thus, the  algebra of formal quantum symmetric function on the
superspaces $\mathbb{R}^{2m|n}$ is isomorphic to
$\sym^{n}(\mathfrak{gl}(\infty)\otimes C(m))$.}
\end{enumerate}
\end{proof}
\subsection{$M$-Weyl algebra}
In this section we introduce the $M$-Weyl algebra. Although
closely related to the Weyl algebra, the $M$-Weyl algebra does not
arises as an instance of the Kontsevich star product.
\begin{defi}
The {\em $M$-Weyl algebra} is the algebra $MW=\mathbb{C} \langle
x,y \rangle[[\hbar]]/ \langle yx-xy- x^{2}\hbar \rangle$. The
letter $M$ stands for meromorphic or mimetic.
\end{defi}
We have the following analogue of Proposition \ref{rweyl}
\begin{prop}\label{rmw}
The map $\rho:MW \longrightarrow \en(\mathbb{C}[x][[\hbar]])$
given by $\rho(x)(f)=x^{-1}f$ and $\rho(y)(f)=-\hbar\frac{\partial
f}{\partial x}$, for any $f\in \mathbb{C}[x][[\hbar]]$ defines a
representation of the $M$-Weyl algebra.
\end{prop}
We order the letters of the $M$-Weyl algebra as follows:
$x<y<\hbar$. Assume we are given $A_i=(a_i,b_i)\in
\mathbb{N}^{2}$, for $i\in[n]$. Set $A=(A_1,\dots,A_n)\in
(\mathbb{N}^{2})^{n}$ and $X^{A_{i}}=x^{a_i}y^{b_i}$, for
$i\in[n]$. Using this notation we have

\begin{defi}
The {\em normal coordinates} $ N_{M}(A,k)$ of \
${\des\prod_{i=1}^{n} X^{A_i}\in MW}$ are defined through the
identity
\begin{equation}\label{mnc}
{\des \prod_{i=1}^{n}X^{A_i}= \sum_{k=0}^{\min}
N_{M}(A,k)x^{|A|+(k,-k)}\hbar^{k}}
\end{equation}
where $0\leq k \leq \min=\min(|a|,|b|)$. For $k>\min$, we set
$N_{M}(A,k)=0$.
\end{defi}

\begin{figure}[h]
\begin{center}
\includegraphics[width=3.5in]{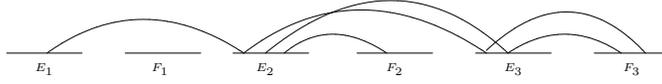}
\caption{Combinatorial interpretation of $N_{M}$.
\label{fig:mero}}
\end{center}
\end{figure}

\begin{thm}\label{emnc}
Let $A,k$  be as in the previous definition, the following
identity holds
\begin{enumerate}
\item{
\begin{equation}\label{cdnor}
N_{M}(A,k)= {\des\sum_{p\vdash k}{b\choose p}
\prod_{i=1}^{n-1}(|a_{>i}|+|p_{>i}|)^{(p_i)}}.
\end{equation}}
\item{Let $E_1,\dots,E_n,F_1,\dots,F_n$ be disjoint sets such that
$\sharp(F_i)=a_i$, $\sharp(E_i)=b_i$, for $i\in[n]$. Set $E=\cup
E_i$, $F=\cup F_i$ and consider the set $M_k$ of all functions
$f:F\longrightarrow P(E)$ such that
\begin{itemize}
\item{$f(x)\cap f(y)=\emptyset$, for all $x,y\in F$.} \item{If
$x\in F_i$, $y\in E_j$, and $y\in f(x)$, then $j<i$.}
\item{${\des\sum_{a\in F}\sharp(f(a))=k}$.}
\end{itemize}
then $N_{M}(A,k)=\sharp(M_k)$.}

\end{enumerate}

\end{thm}
\begin{proof} Using induction one show that the following identity
hold in the $M$-Weyl algebra
\begin{equation}\label{mwal}
y^{b}x^{a}=\sum_{k=0}^{\min}{b\choose k}a^{(k)} x^{a+k}
y^{b-k}\hbar^{k},
\end{equation}
 where $\min=\min(a,b)$ . Several
application of the identity (\ref{mwal}) imply $a.$ Notice that
for given sets $E,F$ such that $|F|=a$ and $|E|=b$, ${b\choose
k}a^{(k)}$ is equal to $\sharp(\{f:F\longrightarrow P(E):\
f(x)\cap f(y)=\emptyset,\
\mbox{for all}\ x,y\in F ,\ \mbox{and}\ \des{\sum_{a\in
F}\sharp(f(a))=k} \})$. This shows $b$, for $n=1$. The general
formula follows from induction. \end{proof} Figure
$\ref{fig:mero}$ illustrate the combinatorial
interpretation of the normal coordinates $N_M(A,6)$ of \\
${\des\prod_{i=1}^{3}X^{A_i}\in MW}$, it shows an example of a
function contributing to $N_M(A,6)$.

\begin{cor} For any given $(t,a,b)\in \mathbb{N}\times\mathbb{N}^{n}\times \mathbb{N}^{n}$, the
following identity holds
$$\prod_{i=1}^{n} (t-|a_{>i}|- |b_{>i}|)_{b_i}=\sum_{p\vdash k}(-1)^{k}{b\choose p}
\prod_{i=1}^{n-1}(|a_{>i}|- |p_{>i}|)^{(p_i)} t_{|b|-k}$$
\end{cor}
\begin{proof} Consider the identity $(\ref{mnc})$ in the representation of
the $M$-Weyl algebra defined in Proposition $\ref{rmw}$. Apply
both sides of the identity $(\ref{mnc})$ to $x^{t}$  and use
Theorem $\ref{emnc}$ formula (\ref{cdnor}).
\end{proof}
The following theorem provides explicit formula for the product of
$m$ elements of the $M$-Weyl algebra. Using the same notation as
in the Theorem $\ref{sp}$, we have the following

\begin{thm}[Symmetric powers of $M$-Weyl algebra]
For any  $A:[m]\!\times\! [n]\longrightarrow \mathbb{N}^{2}$, the
following identity
\begin{equation}
{\des (n!)^{m-1}\overline{\prod_{i=1}^{m}\left(\prod_{j=1}^{n}
X_j^{A_{ij}}\right)}=\sum_{\sigma,k,p}\left(\prod_{i,j} {
b_j^{\sigma} \choose p^{j} }
(|{(a_{j}^{\sigma})}_{>i}|+|p^{j}_{>i}|)^{(p_i^{j})}\right)
\overline{\prod_{j=1}^{n} X_j^{|A_j^{\sigma}|-(k_j,k_j)}}
\hbar^{|k|}}
\end{equation} where
$\sigma\in {\{\id\} \times S_n^{m-1}}$, $k\in\mathbb{N}^{n}$,
$(i,j)\in [m-1]\times [n]$ and $p=p_i^{j}\in
(\mathbb{N}^{m-1})^{n}$, holds in $\sym^{n}(MW)$.
\end{thm}
\begin{proof} By Theorem \ref{pfu} and Theorem \ref{emnc}, we have
\begin{eqnarray*}
{\des (n!)^{m-1}\overline{\prod_{i=1}^{m}\left(\prod_{j=1}^{n}
X_j^{A_{ij}}\right)}}&=&\sum_{\sigma \in \{\id\} \times
S_n^{m-1}}\overline{\prod_{j=1}^{n}\left( \prod_{i=1}^{m}
X_j^{A_{i
\sigma^{-1}_i(j)}}\right)}\\
\mbox{}&=&\sum_{\sigma \in \{\id\} \times
S_n^{m-1}}\overline{\prod_{j=1}^{n}\left(\sum_{k=0}^{\min_j}
N_M(A_j^{\sigma},k)X_j^{|A_j^{\sigma}|+(k,-k)} \hbar^{k}  \right)}\\
\mbox{}&=&\sum_{\sigma,k}\left(\prod_{j=1}^{n}N_M(A_j^{\sigma}
,k_j)  \right)  \overline{\prod_{j=1}^{n}
X_j^{|A_j^{\sigma}|+(k_j,-k_j)}}
\hbar^{|k|}\\
\mbox{} &=&\sum_{\sigma,k,p}\left(\prod_{i,j} { b_j^{\sigma}
\choose p^{j} }
(|{(a_{j}^{\sigma})}_{>i}|+|p^{j}_{>i}|)^{(p_i^{j})}\right)
\overline{\prod_{j=1}^{n} X_j^{|A_j^{\sigma}|-(k_j,k_j)}}
\hbar^{|k|},
\end{eqnarray*}
where $\min_j=\min (|a_j^{\sigma}|,|b_j^{\sigma}|)$.
\end{proof}

\subsection{Cohomological interpretation of symmetric functions}

Let $G$ be a finite group acting on a compact differentiable
manifold $X$.  $G$ acts on $\Ho (X)$, the singular homology groups
with complex coefficients of $X$, as follows:
$$\begin{array}{ccc}
  G\times \Ho (X) & \longrightarrow & \Ho (X) \\
  (g,\alpha) & \longmapsto & g_\ast(\alpha)
\end{array} $$ where $g_\ast(\alpha)$ denotes the push-forward of $\alpha$ by the map $g$.
Similarly $G$ acts on $\cH (X)$, the singular cohomology groups
with complex coefficients of $X$, as follows:
$$\begin{array}{ccc}
  G\times \cH (X) & \longrightarrow & \cH (X) \\
  (g,\alpha) & \longmapsto & g^{\ast}(\alpha)
\end{array} $$ where $g^{\ast}(\alpha)$ denotes the pull-back of $\alpha$ by the map $g$.
It is well-known that $\cH (X/G)=\cH (X)^{G}$, see $\cite{Hir}$.
Identifying $\cH (X)^{G}$ with $\cH (X)_{G}$, we obtain $\cH
(X/G)=\cH (X)_{G}$. Consider $X^{n}/ G^{n}\rtimes K$ where
$K\subset S_n$ and $X^{n}=X\times X\times \dots \times X$. We have
that
$$\cH (X^{n}/G^{n}\rtimes K)\cong \cH (X)^{\otimes n}/G^{n}\rtimes K\cong P_{G,K}(\cH (X)).$$
Next theorem shows how symmetric and supersymmetric functions
arise as the cohomology groups of global orbifolds (quotient of
manifolds by finite group actions). We denote by
$\mathbb{CP}^{\infty}$ the inductive limit of the complex
projective spaces $\mathbb{CP}^{n}$. Also, we let $S^{1}$ be the
unit circle in $\mathbb{C}$.

\begin{thm}\label{last}

\begin{enumerate}
\item{$\cH (((\mathbb{CP}^{\infty})^{m})^{n}/S_n)\cong
\sym_{A_n}(m)$.} \item{$\cH
(((S^{1})^{m})^{n}/S_n)\cong\sym^{n}(\bigwedge[\theta_1,\dots,\theta_m])$.}
\item{$\cH (((\mathbb{CP}^{\infty})^{m}\times
(S^{1})^{k})^{n}/S_n)\cong
\sym^{n}(\mathbb{C}[x_1,\dots,x_n]\otimes \bigwedge
[\theta_1,\dots,\theta_k])$.}
\end{enumerate}
\end{thm}
\begin{proof}
Recall that $\cH(\mathbb{CP}^{n})=\mathbb{C}[x]/(x^{n})$, see
\cite{RB}, which implies that $\cH(\mathbb{CP}^{\infty})=\mathbb{C}[x]$. Thus, by the
remarks above
\begin{eqnarray*}
\cH (((\mathbb{CP}^{\infty})^{m})^{n}/S_n )&\cong&
((\cH(\mathbb{CP}^{\infty})^{\otimes m})^{\otimes
n})_{S_n}\cong\cH
((\mathbb{C}[x]^{\otimes m})^{\otimes n})_{S_n}\\
\mbox{}& \cong & (\mathbb{C}[x_1,\dots,x_m])^{\otimes n}/S_n \cong
\sym_{A_n}(m).
\end{eqnarray*}
Since $\cH(S^{1})\cong\mathbb{C}[\theta]/ \langle
\theta^{2}\rangle$ as graded superalgebras with $\theta$ of
degree $1$, then
$$\cH(((S^{1})^{m})^{n}/S_n)\cong
((\cH(S^{1})^{\otimes m})^{\otimes n})_{S_n}\cong
(\bigwedge[\theta_1,\dots,\theta_m]^{\otimes n})_{S_n}\cong
\sym^{n}(\bigwedge[\theta_1,\dots,\theta_m]).$$ Finally
$$\cH (((\mathbb{CP}^{\infty})^{m}\times
(S^{1})^{k})^{n}/S_n)\cong (\cH(\mathbb{CP}^{\infty})^{\otimes m}
\otimes \cH(S^{1})^{\otimes k})^{\otimes n}/S_n\cong
\sym^{n}(\mathbb{C}[x_1,\dots,x_m]\otimes \bigwedge
[\theta_1,\dots,\theta_k]).$$

\end{proof}

Theorem \ref{last} together  with the quantizations of
$\sym_{A_n}(2m)$ and of
$\sym^{n}(\bigwedge[\theta_1,\dots,\theta_m])$ provided in
Sections \ref{fstw} and \ref{qspf} respectively, give a quantum
product on the cohomology of
$((\mathbb{CP}^{\infty})^{2m})^{n}/S_n$ and $
((S^{1})^{m})^{n}/S_n$ respectively . Notice that the quantum
product on  $\cH (((\mathbb{CP}^{\infty})^{m})^{n}/S_n)$ is
non-commutative and therefore is different to the quantum
cohomology product defined for example in \cite{DMc}.
\subsection*{Acknowledgment}

We thank Nicolas Andruskiewitsch for his advices and
encouragement. We also thank Delia Flores de Chela.

\bibliographystyle{amsplain}
\bibliography{qsf1}

$$\begin{array}{c}
\mbox{Rafael D\'\i az. Instituto Venezolano de Investigaciones Cient\'\i ficas (IVIC).} \ \  \mbox{\texttt{radiaz@ivic.ve}} \\
\!\!\!\!\!\mbox{Eddy Pariguan. Universidad Central de Venezuela (UCV).} \ \  \mbox{\texttt{eddyp@euler.ciens.ucv.ve}} \\
\end{array}$$

\end{document}